\documentclass[11pt,reqno,a4paper]{amsart}
\usepackage{amsmath,amssymb,amsthm,mathtools,mathrsfs}
\usepackage[margin=2.2cm]{geometry}
\usepackage{enumitem}
\usepackage{needspace}
\usepackage[hidelinks]{hyperref}
\numberwithin{equation}{section}
\newtheorem{theorem}{Theorem}[section]
\newtheorem{proposition}[theorem]{Proposition}
\newtheorem{lemma}[theorem]{Lemma}
\newtheorem{corollary}[theorem]{Corollary}
\newtheorem{theoremletter}{Theorem}

\theoremstyle{definition}
\newtheorem{definition}[theorem]{Definition}
\newtheorem{notation}[theorem]{Notation}
\theoremstyle{remark}
\newtheorem{remark}[theorem]{Remark}
\newtheorem{example}[theorem]{Example}

\newcommand{\Fq}{\mathbb F_q}
\newcommand{\Cinf}{\mathbb C_\infty}
\newcommand{\pC}{\widetilde\pi}
\newcommand{\Nq}{\mathbb{N}_q}
\newcommand{\Tr}{\operatorname{Tr}}
\newcommand{\Hcal}{\mathcal H}

\newcommand{\Dcal}{\mathcal D}
\newcommand{\Res}{\operatorname{Res}}
\newcommand{\norm}[1]{\left\lVert#1\right\rVert}
\newcommand{\prodext}{\mathop{\boxtimes}}
\DeclareMathOperator{\ord}{ord}
\usepackage{xcolor}
\newif\ifshowchanges
\showchangestrue

\title[Multivariate Pellarin $L$-values on algebraic curves]{Multivariate Pellarin $L$-values on algebraic curves}
\author{Chuangqiang Hu}
\address{School of Mathematics, Sun Yat-Sen University, Guangzhou, 510275, P.~R.~China}
\email{ \href{mailto:huchq5@mail.sysu.edu.cn}{huchq5@mail.sysu.edu.cn}}

\author{Stephen S.-T. Yau$^*$}
\thanks{ $^*$ Corresponding author.}
\address{Department of Mathematical Sciences, Tsinghua University, Beijing, 100084, P.~R.~China}
\email{ \href{mailto:yau@uic.edu}{yau@uic.edu}}

\author{Lishan Yu}
\address{Beijing Institute of Mathematical Sciences and Applications, Beijing, 101408, P.~R.~China}
\email{ \href{mailto:lishany@bimsa.cn}{lishany@bimsa.cn}}

\begin{document}
\begin{abstract}
We establish a generating trace formula for multicharacter ideal-type $L$-values associated with rank-one Drinfeld modules. For a sign-normalized Hayes module with free period lattice, the formula realizes the $L$-values at weights $n\equiv s\pmod{q^{\deg_\infty}-1}$ as a trace over the ordinary Hilbert class field. We further characterize the polynomials in the exponential variable $X$ arising in multicharacter lattice generating identities over arbitrary function rings. Their coefficients are given explicitly via the shtuka basis, which gives concrete control over their supports and degrees. This completely settles an open problem of Pellarin and Perkins: these polynomials are all constant if and only if $s<2q$. The results on coefficients apply to curves of arbitrary genus, infinite places of arbitrary degree, and nonprincipal period lattices.
\end{abstract}


\maketitle
\enlargethispage{\baselineskip}
\addtocontents{toc}{\protect\setlength{\baselineskip}{12pt}}
\tableofcontents
\section{Introduction}\label{sec:result}
\subsection{Relation to earlier work}

Pellarin's weight-one identity
\cite[Theorem 1]{Pellarin12} expresses a one-variable deformation of
the Carlitz zeta value in terms of the Carlitz period and the special function of Anderson and Thakur
\cite[Section 2.5]{AT90}. Angl\`es and Pellarin subsequently established functional identities for multivariable values \cite{AP14} and connected specializations of the special function and its higher derivatives with Gauss--Thakur sums and abelian extensions \cite{AP15}.

For $1\leqslant s\leqslant2(q-1)$, Perkins obtained explicit generating identities \cite[Theorem 1.1 and Corollary 1.3]{Perkins14}. Pellarin and Perkins then constructed the canonical correction polynomials for arbitrary $s$ \cite[Theorem 2 and Section 4.1]{PP16}.

In the Carlitz case, when $s\geqslant q$ and $s\equiv1\pmod{q-1}$, the constant terms of the corrections recover, up to normalization, the weight-one polynomials $\mathbb B_s$ of Angl\`es, Pellarin, and Tavares Ribeiro
\cite[Section 7.1.1]{APT}. Le and Ngo Dac studied the coefficients of these Bernoulli-type polynomials and deduced Frobenius-twist identities for the weight-one values when $q>(s-1)/(q-1)$
\cite[Theorem 1.3]{LeNgoDac21}.

Goss extended Pellarin's one-variable $L$-series to the
characteristic-$p$ weight space
\cite{Goss13}. Pellarin's tame-series theory further describes the analytic structure and growth of Perkins series
\cite[Sections 4 and 6.4 of the preprint edition]{Pellarin25}. Related polynomial constructions include Perkins's multivariate special polynomials, whose exact degrees are computed in \cite{PerkinsDegree}, and the polynomials interpolating finite twisted harmonic sums of Pellarin and Perkins \cite{PPHarmonic}.

Over a general function ring, the ideal-series approach originates in Anderson's work on $A$-harmonic series and log-algebraicity \cite{Anderson94,Anderson96}. Angl\`es, Ngo Dac, and Tavares Ribeiro developed several-variable twisted zeta functions in arbitrary genus and for infinite places of arbitrary degree, and proved their analytic continuation
\cite[Section 4]{ANDTR16}. They subsequently related equivariant $A$-harmonic series to Stark units and established several-variable log-algebraicity
\cite[Sections 3 and 4]{ANDTRStark17}. Their construction of special functions from shtuka functions also yields a normalized multivariable weight-one rationality theorem for
$s\equiv1\pmod{q^{\deg\infty}-1}$ \cite[Theorem 5.6]{ANDTR17}. When the infinite place is rational, Green and Papanikolas obtained explicit weight-one element and ideal formulas on elliptic curves
\cite[Theorems 7.1 and 7.3]{GP18}. Also at a rational infinite place, Ferraro established functional identities for ideal lattices and an Anderson zeta trace identity in arbitrary genus
\cite[Theorems 7.27 and 8.9]{Ferraro24}.

Poonen established a duality between the kernels of additive series and their adjoints \cite{Poonen96}. El-Guindy and Papanikolas derived explicit identities for Anderson generating functions over polynomial rings \cite{EP14}. Gazda and Maurischat constructed special functions for Anderson $A$-modules \cite{GM21} and later recovered periods from residues at the characteristic graph \cite{GM24}. Ferraro proved the rationality of the pairing between universal special and dual special functions for Drinfeld modules of arbitrary rank
\cite[Theorem 7.3]{Ferraro25}. Hu, Yau, and Yu constructed the normalized Cauchy realization of the rank-one generating differential
\cite[Definition 4.19 and Theorem 4.20]{HYYCauchy}.
We use its lattice identity, Frobenius equation, and residue
pairing in the form developed in \cite[Sections 3.1--3.4 and Theorem 4.6]{HYY}. Pellarin's recent account
\cite{Pellarin26} places these special-function identities within the broader theory of zeta functions over curves.

The multicharacter values studied here are the negatives of the
Taylor coefficients of the lattice generating series. For $s=1$, they
are the higher-weight differentials of \cite{HYY}; for the Carlitz module,
they recover the multicharacter values of Pellarin--Perkins.

\subsection{The motivating question}
\leavevmode
Over $\Fq[\theta]$, Pellarin and Perkins expressed the generating
function for multicharacter $L$-values in terms of Anderson generating
functions and canonical correction polynomials. The characters are the
evaluation maps $\chi_{t_i}(a)=a(t_i)$ at independent variables
$t_1,\ldots,t_s$. Their Question 3 asks whether all these corrections
are constant in the exponential variable \cite{PP16}.

We determine every correction coefficient for every $s\geqslant1$ and
prove that all the corrections are constant if and only if $s<2q$.
For a nonempty factor set $S\subseteq\{1,\ldots,s\}$, write
$\mathfrak g_{S,d}$ for the coefficient of $X^d$ in the differential
normalization introduced below. Our formula is
\[
 \mathfrak{g}_{S,d}
 =\sum_{(n_i)\in\mathbb{E}_{S,d}}
   \prodext_{i\in S}\mathfrak{B}_i^{[n_i]}
 \qquad(d\geqslant0).
\]
Here $\mathfrak B^{[n]}$ is the basis differential constructed in
Section~\ref{sec:basis}, and $\mathbb E_{S,d}$ consists of the tuples of
nonnegative integers satisfying $\sum_{i\in S}q^{-n_i}=q(d+1)$.
Theorem~\ref{thm:formula} shows that this sum is finite and that its
terms are distinct elements of an external tensor basis. Consequently,
$\mathfrak g_{S,d}\ne0$ if and only if
\[
 |S|=q(d+1)+(q-1)k
 \quad\text{for some integer }k\geqslant0.
\]

\subsection{Multicharacter values and class-field trace identities}

We work over a smooth, projective, geometrically connected curve
${\mathcal X}/\Fq$. Fix a closed point $\infty$ of degree $N\geqslant1$,
and set
$A:=\Gamma({\mathcal X}\setminus\{\infty\},\mathcal O_{\mathcal X})$ and
$K:=\Fq({\mathcal X})$. The completion of $K$ at $\infty$ is denoted by
$K_\infty$, with residue field $\mathbb F_\infty:=\kappa(\infty)$, and
$\Cinf$ denotes the completion of an algebraic closure of $K_\infty$.
Throughout, we fix an isometric $\Fq$-embedding
$\iota:K_\infty\hookrightarrow\Cinf$, which serves as the coefficient
embedding.

Consider a generic-characteristic rank-one Drinfeld $A$-module $\varphi$
over $\Cinf$, with structural map $\iota|_A$, period lattice $\Lambda$, and
exponential $e(w)$. For a fixed integer $s\geqslant1$, take $s$ independent
copies of ${\mathcal X}$. On the $i$-th factor, $\Hcal_i$ denotes the entire
generating differential of \cite[Section 3.1 and Theorem 4.6]{HYY}, while
$\prodext$ denotes the external tensor product.

The multicharacter Pellarin--Perkins generating function is defined by
\[
 \mathscr C_\varphi^s(w)
 :=\sum_{0\ne\lambda\in\Lambda}
       \frac{\prodext_{i=1}^s\Hcal_i(\lambda)}{w-\lambda},
\]
and the associated multicharacter Ferraro series of weight $n\geqslant1$,
in differential normalization, is
\[
C_{\varphi,n}^{s}:=\sum_{0\ne\lambda\in\Lambda}
       \frac{\prodext_{i=1}^s\Hcal_i(\lambda)}{\lambda^n}.
\]

\Needspace{6\baselineskip}
\begin{theoremletter}
\label{thm:introtrace}
For every rank-one $\varphi$, the expansion
\[
 \mathscr C_\varphi^s(w)
 =-\sum_{\substack{n\geqslant1\\n\equiv s\ (q-1)}}C_{\varphi,n}^{s}w^{n-1}
\]
holds when $|w|<\min_{0\ne\lambda\in\Lambda}|\lambda|$.
Its coefficients are rational external differentials over the field of
definition of $\varphi$ and vanish unless $n\equiv s\pmod{q-1}$.
\end{theoremletter}

For the ideal-type $L$-value statement, we specialize to a
sign-normalized Hayes module $\psi$ with free period lattice. Choosing a
generator $\pC\in\Cinf^\times$ gives $\Lambda_\psi=\pC\iota(A)$. We view
the sign and ordinary Hilbert class fields $H_A^+$ and $H_A$ inside $\Cinf$
through compatible embeddings extending $\iota|_K$.
To describe the ideal series, let $\mathfrak a$ be a nonzero integral ideal
and put $d_{\mathfrak a}:=\partial\psi_{\mathfrak a}$ and
$u_{\mathfrak a,i}:=p_i^*(\psi_{\mathfrak a}(\omega_\psi)/\omega_\psi)$.
Here $\psi_{\mathfrak a}:\psi\to\mathfrak a*\psi$ is the monic ideal isogeny
for the Artin action \cite[Section 3.5]{HYY}, $\partial$ extracts the
$\tau^0$-coefficient, $p_i$ is the $i$-th projection, and $\omega_\psi$ is
a nonzero special function as in \cite[Section 2.1]{HYY}. With these
conventions, define the associated ideal-type $L$-value and period product by
\[
 \mathbb L_{\psi}^s(n):=\sum_{0\ne\mathfrak a\subset A}
       \frac{\prod_{i=1}^s u_{\mathfrak a,i}}{d_{\mathfrak a}^{\,n}},
 \qquad
 \mathscr W_{\psi,s}:=\prodext_{i=1}^s\Hcal_{\psi,i}(\pC).
\]
For each $n\geqslant1$, define the companion element series by
\[
 {L^s}(A;n)
 :=\sum_{0\ne a\in A}\frac{\prod_{i=1}^s\chi_i(a)}{\iota(a)^n},
\]
where $\chi_i(a)$ denotes the geometric copy of $a$ on the $i$-th factor.
The multipliers $u_{\mathfrak a,i}$ are independent of the choice
of $\omega_\psi$. Proposition~\ref{prop:mc-class-decomposition}(i)
proves convergence of the ideal series for every $n\geqslant1$.

\Needspace{13\baselineskip}
\begin{theoremletter}
\label{thm:introidealtrace}
In this Hayes setting:
\begin{enumerate}[label=(\roman*)]
\item For every $n\geqslant1$,
\[
 C_{\psi,n}^{s}=\pC^{-n}\mathscr W_{\psi,s}\,
 {L^s}(A;n).
\]
Near $z=0$, the corresponding generating expansion is
\begin{equation}\label{eq:intro-element-generating}
 \mathscr C_\psi^s(\pC z)
 =-\frac{\mathscr W_{\psi,s}}{\pC}
 \sum_{\substack{n\geqslant1\\n\equiv s\ (q-1)}}
 {L^s}(A;n)z^{n-1}.
\end{equation}
\item Identify $\mathbb F_\infty$ with its constant-field lift in
$K_\infty$. For a power series $B(w)$, define
\[
 (\mathcal P_{q^N}^{s}B)(w):=\frac1{q^N-1}
     \sum_{\varepsilon\in\iota(\mathbb F_\infty^\times)}
          \varepsilon^{1-s}B(\varepsilon w).
\]
Then, near $z=0$,
\begin{equation}\label{eq:intro-mc-trace}
   \left.\Tr_{H_A/K}
      \bigl((\mathcal P_{q^N}^{s}\mathscr C_{\psi}^s)(w)\bigr)
   \right|_{w=\pC z} = \frac{\mathscr W_{\psi,s}}{\pC}\sum_{\substack{n\geqslant1\\n\equiv s\ (q^N-1)}}\mathbb L_{\psi}^s(n)z^{n-1}.
\end{equation}
The trace is taken on the Taylor coefficients before the period
substitution, with $w$ and the geometric variables held fixed. The resulting
finite sum over ideal classes also provides a meromorphic continuation of the
series on the right.
\end{enumerate}
\end{theoremletter}

When $s=1$, Theorem~\ref{thm:introtrace} and
Theorem~\ref{thm:introidealtrace}(ii) reduce, respectively, to
\cite[Theorems D and E(i)]{HYY}.

\subsection{Canonical corrections over an arbitrary curve}

\begin{definition}\label{def:correction-polynomials}
For a nonempty $S\subseteq\Sigma_s:=\{1,\ldots,s\}$, consider a polynomial
$\mathfrak G_S(X)=\sum_{d}\mathfrak{g}_{S,d}X^d$ whose coefficients are
rational external differentials on the factors indexed by $S$. We call
it a \emph{canonical correction polynomial} if
\[
 \mathfrak G_S(e(w))=\sum_{\lambda\in\Lambda}
       \frac{\prodext_{i\in S}\Hcal_i(\lambda-w)}{w-\lambda}
\]
holds for every $w\notin\Lambda$.
\end{definition}

\begin{theoremletter}
\label{thm:intromain}
\textup{(i)} For every nonempty $S\subseteq\Sigma_s$, there exists a unique
canonical correction polynomial $\mathfrak G_S(X)$ in the sense of
Definition~\ref{def:correction-polynomials}.
\par\medskip\noindent\textup{(ii)} With $\mathfrak{G}_\varnothing(X):=X^{-1}$, the generating
function satisfies the following meromorphic identity in $w$:
\[
 \mathscr C_\varphi^s(w) 
 =\sum_{S\subseteq\Sigma_s}
       \mathfrak{G}_S(e(w))\boxtimes\prodext_{i\notin S}\Hcal_i(w).
\]
\par\medskip\noindent\textup{(iii)} Regarding the zero polynomial as constant, $\mathfrak{G}_S$ is
constant if and only if
\begin{equation}\label{intro:constancy}
  |S|\notin\Nq :=  
   q\mathbb Z_{\geqslant2}+(q-1)\mathbb Z_{\geqslant0}.
\end{equation}
In particular, $\mathfrak{G}_S$ is constant for every nonempty
$S\subseteq\Sigma_s$ if and only if $s<2q$.
\end{theoremletter}

\subsection{Organization of the paper}

\leavevmode
Section~\ref{sec:preliminaries} fixes the period and
generating-differential normalizations and constructs the shtuka basis.
Section~\ref{sec:geometry} constructs the correction polynomials by finite
interpolation and passage to the period lattice.
Section~\ref{sec:coeff-formula} derives the finite coefficient formula,
determines the exact support and degrees, and descends the coefficients
to the module's definition field. These two sections prove
Theorem~\ref{thm:intromain}. Section~\ref{sec:mc-traces} proves the
rationality and covariance of the multicharacter values and the
Hilbert-class trace identities in Theorems~\ref{thm:introtrace}
and~\ref{thm:introidealtrace}. Section~\ref{sec:curve-examples} gives
Carlitz and nonprincipal-lattice examples, together with formulas for
repeated parameters.

\section{Preliminaries on period lattices and generating differentials}\label{sec:preliminaries}

Retain the curve ${\mathcal X}$, ring $A$, and field $\Cinf$ from the
introduction. Write $q:=p^e$ with $p$ prime, put
$U:={\mathcal X}\setminus\{\infty\}$, and denote the absolute value on
$\Cinf$ by $|\cdot|$. For nonempty $S\subseteq\Sigma_s$, write
$\mathcal X^S=\prod_{i\in S}\mathcal X_i$, where $\mathcal X_i$ is
the $i$-th independent copy of $\mathcal X$.

\subsection{Period and shtuka data}

For the background on Drinfeld modules and sign-normalized
rank-one modules, we refer to \cite{Hayes92} and
\cite[Chapters 4 and 7]{Goss96}.
\begin{definition}
With the coefficient embedding $\iota$ fixed above, a rank-one Drinfeld
$A$-module of generic characteristic is an $\Fq$-algebra homomorphism
\[
      \varphi : A \to \Cinf \{ \tau \},
\]
 such that
$\partial\varphi_a=\iota(a)$ for every $a\in A$ and
$\deg_\tau\varphi_a=\dim_{\Fq}(A/aA)$ for every $a\ne0$.
Here $\tau$ is the Frobenius variable, with
$\tau c=c^q\tau$, and $\partial$ extracts the $\tau^0$-coefficient.
\end{definition}
\begin{notation}\label{not:lattice-data}
Write $e(w):=\exp_\varphi(w)$ for the exponential function
associated with $\varphi$, so that $e'(w)=1$, and let $\Lambda$ be its
period lattice.
\end{notation}

The symbols $(j)$ denote coefficient Frobenius: constants are
raised to $q^j$, with the induced action on geometric points and divisors,
while functions and differentials defined over $\Fq$ are fixed. Inverse
twists exist over $\Cinf$.

For the shtuka description of rank-one Drinfeld modules,
see \cite[Section 0.3]{Thakur93} and
\cite[Section 2.3]{ANDTR17}.
\begin{notation}\label{not:shtuka-data}
Write $g$ for the genus of ${\mathcal X}$.
The standard geometric data consist of an effective divisor $V_\varphi$,
a geometric point $P$ above $\infty$, and the shtuka function $f$, with 
\begin{equation}\label{eq:geometric-data}
 \begin{gathered}
 \deg V_\varphi=g,\quad h^0(V_\varphi)=1,\quad
 \Xi,P^{(-1)}\notin\operatorname{supp}V_\varphi,\\
 (f)=V_\varphi^{(1)}-V_\varphi+\Xi-P.
 \end{gathered}
\end{equation}
Here $\Xi$ is the characteristic point. The normalized
shtuka differential $\Delta_\varphi$ is the unique rational differential on
${\mathcal X}_{\Cinf}$ satisfying
\begin{equation}\label{eq:normalizations}
 (\Delta_\varphi)\geqslant V_\varphi-P^{(-2)}-P^{(-1)},\qquad
 \Res_\Xi\!\left(\frac{\Delta_\varphi^{(1)}}{f}\right)=1.
\end{equation}
Its existence and uniqueness follow from
\cite[Definition 6.1 and Theorem 6.2]{HYYCauchy}. The associated Cauchy
differential is $C_\varphi:=\Delta_\varphi^{(1)}/f$.
\end{notation}

\begin{proposition}\label{prop:shtuka-log-cauchy}
\textup{(i)} If $\log_\varphi$ denotes the formal composition inverse of
$e=\exp_\varphi$, then
\begin{equation}\label{eq:shtuka-log}
 \log_\varphi(Z)=\sum_{i\geqslant0}\Res_\Xi\!\left(
 \frac{\Delta_\varphi^{(i+1)}}{f f^{(1)}\cdots f^{(i)}}\right)Z^{q^i},
\end{equation}
where the denominator is $f$ for $i=0$.
\par\medskip\noindent\textup{(ii)} The Cauchy differential
$C_\varphi$ satisfies $(C_\varphi)\geqslant V_\varphi-\Xi-P^{(-1)}$.
It has simple poles at $\Xi$ and $P^{(-1)}$ with residues $1$ and $-1$,
respectively.
\end{proposition}

\begin{proof}

\emph{(i).} This is the residue formula for the formal logarithm with
the normalized shtuka differential
\cite[Definition 6.24 and Theorem 6.26]{HYYCauchy}.

\emph{(ii).} Twisting the divisor bound in
\eqref{eq:normalizations} and subtracting $(f)$ gives
$(C_\varphi)\geqslant V_\varphi-\Xi-P^{(-1)}$. The residue normalization in
\eqref{eq:normalizations} gives residue one at $\Xi$, and the divisor bound
leaves no other possible pole except $P^{(-1)}$. The global residue theorem
forces a simple pole of residue $-1$ there.
\end{proof}

\subsection{Generating differentials and uniform bounds}

\begin{notation}\label{not:generating-differential}
Put $\Omega_A:=\Gamma(U,\Omega^1_{{\mathcal X}/\Fq})$.
Both $\Lambda$ and $\Omega_A$ are projective $A$-modules of rank one.
Write $\mathbb M_A:=\Cinf\widehat\otimes_{\Fq}\Omega_A$
for the Gauss completion of $\Cinf\otimes_{\Fq}\Omega_A$,
and $\chi(a):=1\otimes a$ for the geometric copy of $a\in A$.
Write $\norm{\cdot}_{\mathbb M_A}$ for the Gauss norm on
$\mathbb M_A$.

For the independent copies introduced above, pullbacks under
the projections ${\mathcal X}^S\to {\mathcal X}_i$ ($i\in S$) are indicated by subscript $i$. Products of differentials on distinct factors
mean external tensor products, over the coefficient field unless indicated
otherwise. The binary symbol $\boxtimes$ and the symbol
$\prodext_{i\in S}$ denote external tensor products on the indicated
factors.
\end{notation}

The entire $\Fq$-linear generating differential $\Hcal(w)$,
constructed in \cite[Definition 4.19]{HYYCauchy} and used in
\cite[Section 3.1 and Theorem 4.6]{HYY}, is normalized by its Cauchy
kernel. It takes values in $\mathbb M_A$.

\begin{proposition}\label{prop:generating-differential-identities}
The generating differential $\Hcal$ satisfies
\begin{equation}\label{eq:inputs}
 \frac{\Hcal(w)}{e(w)}
 =\sum_{0\ne\lambda\in\Lambda}\frac{\Hcal(\lambda)}{w-\lambda},
 \qquad
 \Hcal(\lambda)^{(1)}=f\Hcal(\lambda),
 \qquad
 \Hcal'(0)=C_\varphi.
\end{equation}
\end{proposition}

\begin{proof}
These are the lattice identity, the Frobenius
eigenvalue property, and the Cauchy normalization of the generating
differential \cite[Section 3.1 and Theorem 4.6]{HYY}.
\end{proof}

Choose $0\ne\pi\in\Lambda$. Since $\Lambda$ has rank one
over $\iota(A)$, the quotient $\pi^{-1}\Lambda$ lies in $\iota(K)$.
The fractional ideal
$\mathfrak j:=\iota^{-1}(\pi^{-1}\Lambda)\subset K$ contains $A$.
If $\Lambda$ is free and $\pi$ is a generator, then $\mathfrak j=A$.
We extend $\chi$ and its pullbacks $\chi_i$ to $K$ by localization.
Identities involving fractional elements are read after localizing the
differential modules at the nonzero geometric functions $\chi_i(b)$,
$b\in A$.
The period and residue-duality identities of
\cite[Sections 3.1--3.4 and Theorem 4.6]{HYY} give the following
period-linearity properties.
\begin{proposition}\label{prop:fractional-linearity}
\textup{(i)} The restriction of $\Hcal$ to $\Lambda$ is $A$-linear:
\[
 \Hcal(\iota(a)\lambda)=\chi(a)\Hcal(\lambda)
       \qquad(a\in A,\ \lambda\in\Lambda).
\]
\par\medskip\noindent\textup{(ii)} For every $a\in\mathfrak j$ and every $i$,
\begin{equation}\label{eq:fractional-period-linearity}
 \Hcal_i(\pi\iota(a))=\chi_i(a)\Hcal_i(\pi).
\end{equation}
\end{proposition}

\begin{proof}
The period map is the inverse of the $A$-linear characteristic-residue
isomorphism of \cite[Section 3.1]{HYY}, which proves (i).
For (ii), choose $0\ne b\in A$ with $ba\in A$. Both $\pi$ and
$\pi\iota(a)$ belong to $\Lambda$, so (i) gives
\[
 \chi_i(b)\Hcal_i(\pi\iota(a))
 =\Hcal_i(\pi\iota(ba))
 =\chi_i(ba)\Hcal_i(\pi).
\]
Division by the nonzero geometric function $\chi_i(b)$ gives
\eqref{eq:fractional-period-linearity}.
\end{proof}

\begin{notation}\label{not:scalar-trivialization}
Fix $0\ne\nu\in\Omega^1_{K/\Fq}$ and put
$\nu_S:=\prodext_{i\in S}\nu_i$. Define
\[
 h_i(w):=\Hcal_i(w)/\nu_i,\qquad c_i:=h_i'(0)=C_{\varphi,i}/\nu_i.
\]
\end{notation}

Let $(F,\norm{\cdot}_F)$ be a complete nonarchimedean
valued extension of $\Cinf$. In the geometric application, $F$ is the
Gauss completion of the fraction field of $\Cinf\otimes_{\Fq}A^{\otimes s}$,
and $\norm{\cdot}_F$ is its Gauss norm.

The geometric Gauss norm is multiplicative because ${\mathcal X}^{\Sigma_s}$ is
geometrically integral. The reductions of two finite sums of norm one
are nonzero in the domain
$\overline{\Fq}\otimes_{\Fq}A^{\otimes s}$, so their product also has
norm one. Scaling, then passing to fractions and completion, proves
multiplicativity on $F$.

Division by $\nu_S$ also defines the scalar realization of completed
differentials. Put
\[
 \mathbb M_S=\Cinf\widehat\otimes_{\Fq}
                \left(\bigotimes_{\Fq,\,i\in S}\Omega_A\right),
\]
the completed external tensor product over $\Cinf$ of the one-factor
modules $\mathbb M_A$. Over $R=A^{\otimes_{\Fq}S}$, division by $\nu_S$
identifies $\bigotimes_{\Fq,\,i\in S}\Omega_A$ with an invertible
fractional ideal $I$. Choose $a_\ell\in I$ and $b_\ell\in I^{-1}$ such that
$\sum_\ell a_\ell b_\ell=1$. After extension of scalars to $\Cinf$, the maps
\[
 x\longmapsto(b_\ell x)_\ell,\qquad
 (y_\ell)_\ell\longmapsto\sum_\ell a_\ell y_\ell
\]
give a split presentation of $\Cinf\otimes_{\Fq}I$.

Every nonzero $a_\ell,b_\ell$ is a geometric rational function of Gauss
norm one, so both maps have norm at most one. The split presentation
computes the module norm of $x$ as
$\max_\ell\norm{b_\ell x}_F=\norm{x}_F$, which is its scalar Gauss norm.
Passing to completions
gives an isometric embedding $\mathbb M_S\to F$, with closed image,
induced by division by $\nu_S$. In particular, the scalarizations
$h_i$ are entire $F$-valued functions.

\Needspace{7\baselineskip}
\begin{lemma}\label{lem:period-bound}
For each $i$, there exists a constant $M_i>0$, independent of $\lambda$,
such that the following bounds hold for every $\lambda\in\Lambda$:
\begin{enumerate}[label=(\roman*)]
\item The period differentials satisfy $\norm{\Hcal_i(\lambda)}_{\mathbb M_A}\leqslant M_i$.
\item The scalar period values satisfy $\norm{h_i(\lambda)}_F\leqslant M_i$.
\end{enumerate}
\end{lemma}

\begin{proof}
Choose $A$-module generators $\lambda_1,\ldots,\lambda_r$ of $\Lambda$.
For $\lambda=\sum_j\iota(a_j)\lambda_j$,
Proposition~\ref{prop:fractional-linearity}(i) gives
\[
 \Hcal_i(\lambda)=\sum_j\chi_i(a_j)\Hcal_i(\lambda_j),\qquad
 h_i(\lambda)=\sum_j\chi_i(a_j)h_i(\lambda_j).
\]
Each geometric coefficient $\chi_i(a_j)$ has Gauss norm at
most one. The nonarchimedean
triangle inequality proves both bounds with
\[
 M_i=\max_{1\leqslant j\leqslant r}
 \left\{\norm{\Hcal_i(\lambda_j)}_{\mathbb M_A},
              \norm{h_i(\lambda_j)}_F\right\}+1.
\]
\end{proof}

\subsection{A differential basis from the shtuka function}
\label{sec:basis}

\begin{samepage}
\begin{lemma}\label{lem:W}
\textup{(i)} The divisor
\begin{equation}\label{eq:W}
 W=(\Delta_\varphi)-V_\varphi+P^{(-2)}+P^{(-1)}
\end{equation}
is effective of degree $g$, satisfies $h^0(W)=1$, and has neither
$\Xi^{(-1)}$ nor $P^{(-2)}$ in its support.
\par\medskip\noindent\textup{(ii)} For $\gamma=(f\Delta_\varphi)^{(-1)}/\Delta_\varphi$,
\begin{equation}\label{eq:gamma-divisor}
 (\gamma)=W^{(-1)}-W+\Xi^{(-1)}-P^{(-3)}.
\end{equation}

\end{lemma}
\end{samepage}

\begin{proof}

\emph{(i).} The divisor $W$ is effective by \eqref{eq:normalizations}, and its degree is
$(2g-2)-g+2=g$. Serre duality gives
\[
 h^1(W)=h^0(V_\varphi-P^{(-2)}-P^{(-1)})=0:
\]
a section of $\mathcal O(V_\varphi-P^{(-2)}-P^{(-1)})$ lies in
$H^0(V_\varphi)=\Cinf$ and vanishes at $P^{(-1)}$, since
$P^{(-1)}\notin\operatorname{supp}V_\varphi$.
Riemann--Roch now gives $h^0(W)=1$.

Proposition~\ref{prop:shtuka-log-cauchy}(ii) gives the
divisor bound and the two residues of $C_\varphi$. Moreover,
\[
 W^{(1)}=(C_\varphi)-V_\varphi+\Xi+P^{(-1)}.
\]
Since $C_\varphi$ has simple poles at $\Xi$ and $P^{(-1)}$ and
$V_\varphi$ avoids those points, $W^{(1)}$ has coefficient zero at
$\Xi$ and $P^{(-1)}$. Twisting back shows that $W$ avoids
$\Xi^{(-1)}$ and $P^{(-2)}$.

\emph{(ii).} Substituting
\eqref{eq:W} and \eqref{eq:geometric-data} into the divisor of
$(f\Delta_\varphi)^{(-1)}/\Delta_\varphi$ gives \eqref{eq:gamma-divisor}.
\end{proof}

Define rational differentials
\begin{equation}\label{eq:B-definition}
 \mathfrak{B}^{[0]}=\Delta_\varphi,\qquad
 \mathfrak{B}^{[n+1]}=f^{(-1)}\left(\mathfrak{B}^{[n]}\right)^{(-1)}\quad(n\geqslant0).
\end{equation}
In particular, $\mathfrak{B}^{[1]}=(f\Delta_\varphi)^{(-1)}$.

\begin{notation}\label{not:Gamma}
Set
\[
 \Gamma^{[0]}:=1,\qquad
 \Gamma^{[n]}:=\prod_{j=0}^{n-1}\gamma^{(-j)}\quad(n\geqslant1).
\]
These functions will be used in Section~\ref{sec:formula} to expand
the normalized correction coefficients.
\end{notation}

\begin{proposition}\label{prop:basis}
The differentials in \eqref{eq:B-definition} form a basis
\begin{equation}\label{eq:differential-basis}
 \Gamma(U_{\Cinf},\Omega^1(-V_\varphi))
       =\bigoplus_{n\geqslant0}\Cinf \mathfrak{B}^{[n]}.
\end{equation}
\end{proposition}

\begin{proof}
Telescoping gives
\begin{equation}\label{eq:Gamma-divisor}
 (\Gamma^{[n]})=W^{(-n)}-W+\sum_{j=1}^n\Xi^{(-j)}-D_n,
 \qquad D_n=\sum_{j=3}^{n+2}P^{(-j)},\quad D_0=0.
\end{equation}
Since $D_n$ is effective, $h^1(W+D_n)=0$ follows from
$h^1(W)=0$. Riemann--Roch then gives $h^0(W+D_n)=n+1$.

For $n\geqslant1$, formula~\eqref{eq:Gamma-divisor} places $\Gamma^{[n]}$
in $H^0(W+D_n)$. All points $\Xi^{(-j)}$ lie in $U$, and
$\operatorname{mult}_{P^{(-n-2)}}W^{(-n)}
=\operatorname{mult}_{P^{(-2)}}W=0$ by Lemma~\ref{lem:W}(i). Since
$D_n=D_{n-1}+P^{(-n-2)}$, we obtain
\[
 \ord_{P^{(-n-2)}}\Gamma^{[n]}
 +\operatorname{mult}_{P^{(-n-2)}}(W+D_{n-1})=-1.
\]
Thus $\Gamma^{[n]}\notin H^0(W+D_{n-1})$. Induction proves that
$\Gamma^{[0]},\ldots,\Gamma^{[n]}$ is a basis of $H^0(W+D_n)$.

The Frobenius orbit of $P$ consists of all geometric points over
$\infty$, and the multiplicity of each such point in $D_n$ tends to
infinity. Every section of $\mathcal O(W)$ on $U_{\Cinf}$ therefore
belongs to $H^0(W+D_n)$ for some $n$. Hence
\[
 \Gamma(U_{\Cinf},\mathcal O(W))
       =\bigoplus_{n\geqslant0}\Cinf\Gamma^{[n]}.
\]
Multiplication by $\Delta_\varphi$ identifies this space with
$\Gamma(U_{\Cinf},\Omega^1(-V_\varphi))$ by \eqref{eq:W}.
Finally $\Delta_\varphi\Gamma^{[n]}=\mathfrak{B}^{[n]}$, as follows from the recurrence
\eqref{eq:B-definition} and the definition of $\gamma$.
\end{proof}

\section{Canonical corrections by lattice interpolation}\label{sec:geometry}

\subsection{Finite linear interpolation and trace bounds}
\label{sec:finite}

The finite interpolation method follows the approach to
products of evaluation characters in \cite[Section 2.4]{Perkins14}.
Here the interpolation space is an arbitrary finite $\Fq$-subspace
of $\Cinf$, and the maps take values in a complete valued field.
\begin{notation}\label{not:finite-interpolation}
Let $V\subset\Cinf$ be a finite $\Fq$-vector space of dimension $D\geqslant1$.
\begin{enumerate}[label=(\roman*)]
\item Define
\begin{equation}\label{eq:EV}
 E_V(w):=w\prod_{0\ne\lambda\in V}(1-w/\lambda)
       =\sum_{j=0}^D \alpha_j w^{q^j}, 
\end{equation}
for some constants $\alpha_j\in\Cinf$.
In particular, $\alpha_0=1$ and
$\alpha_D=\alpha\ne0$.
\item For any $\Fq$-linear maps
$h_i:V\to F$, define their interpolation polynomials by
\begin{equation}\label{eq:IV}
 I_{V,i}(w):=E_V(w)\sum_{0\ne\lambda\in V}
                       \frac{h_i(\lambda)}{w-\lambda}
       =\sum_{j=0}^{D-1}b_{V,i,j}w^{q^j}.
\end{equation}
\end{enumerate}
\end{notation}

\begin{lemma}\label{lem:interpolation-values}
The polynomial $I_{V,i}$ defined by the first expression in
\eqref{eq:IV} has degree less than $q^D$ and satisfies
$I_{V,i}(\lambda)=h_i(\lambda)$ for every $\lambda\in V$.
\end{lemma}
\begin{proof}
Each summand is a polynomial of degree $q^D-1$, so $I_{V,i}$ has
degree less than $q^D$. Let $\lambda_0\in V$ be nonzero. As
$w\to\lambda_0$, the summand with $\lambda=\lambda_0$ tends to
$h_i(\lambda_0)E_V'(\lambda_0)=h_i(\lambda_0)$, while every other
summand vanishes. Thus $I_{V,i}(\lambda_0)=h_i(\lambda_0)$. At $w=0$, every
summand vanishes, and $I_{V,i}(0)=0=h_i(0)$ by $\Fq$-linearity.
\end{proof}
To obtain the second expression in \eqref{eq:IV}, choose an
$\Fq$-basis $v_1,\ldots,v_D$ of $V$. The Moore matrix
$(v_a^{q^j})_{1\leqslant a\leqslant D,\,0\leqslant j<D}$ is invertible:
a nonzero vector in its kernel would give a nonzero $q$-polynomial
of degree at most $q^{D-1}$ vanishing at all $q^D$ points of $V$.
The resulting interpolant agrees with $h_i$ on the basis, hence on
$V$ by $\Fq$-linearity. Lemma~\ref{lem:interpolation-values} and
uniqueness of interpolation in degree less than $q^D$ identify it with $I_{V,i}$.

\begin{notation}
Throughout, $S$ denotes a nonempty subset of
$\Sigma_s$, and $m:=|S|$.
\end{notation}

\begin{lemma}\label{lem:finite}
\textup{(i)} For $m:=|S|\geqslant1$, the finite sum
\[
 \sum_{\lambda\in V}
     \frac{\prod_{i\in S}I_{V,i}(\lambda-w)}{w-\lambda}
\]
is $G_{V,S}(E_V(w))$ for a polynomial $G_{V,S}\in F[X]$
of degree at most $\lfloor m/q\rfloor-1$, with a negative bound
meaning that the polynomial is zero.
\par\medskip\noindent\textup{(ii)}  If the coefficient of $X^d$ in $G_{V,S}$ is nonzero,
then $m=q(d+1)+(q-1)k$ for some integer $k\geqslant0$.

\end{lemma}

\begin{proof}

\emph{(i).} In the finite free \'etale $F[X]$-algebra
$B:=F[X,Y]/(E_V(Y)+X)$, set
\begin{equation}\label{eq:trace}
 G_{V,S}(X):=-\Tr_{B/F[X]}
            \left(\frac{\prod_{i\in S}I_{V,i}(Y)}{Y}\right).
\end{equation}
The argument of the trace is a polynomial because $S$ is nonempty.
At $X=E_V(w)$ the roots of $E_V(Y)+X$ are $\lambda-w$, $\lambda\in V$,
which proves the finite sum identity.

Put $t_n(X):=\Tr_{B/F[X]}(Y^n)$.
To compute these traces, put $P_X(Y):=\alpha^{-1}(E_V(Y)+X)$.
Over a splitting field of $F(X)$, the trace is the sum over the
distinct roots $y$ of $P_X$. Hence, as formal series at $T=0$,
\[
 \sum_{n\geqslant0}t_n(X)T^n
 =\sum_{P_X(y)=0}\frac1{1-y T}
 =\frac{T^{-1}P_X'(T^{-1})}{P_X(T^{-1})}.
\]
Since $P_X'=\alpha^{-1}$, multiplying numerator and denominator by $T^{q^D}$
gives
\begin{equation}\label{eq:resolvent}
 \sum_{n\geqslant0}t_n(X)T^n
 =\frac{\alpha^{-1}T^{q^D-1}}
  {1+\displaystyle\sum_{j=0}^{D-1}(\alpha_j/\alpha)T^{q^D-q^j}
       +(X/\alpha)T^{q^D}}.
\end{equation}
Consequently $t_n=0$ for $n<q^D-1$ and otherwise
\[
 \deg_Xt_n\leqslant\left\lfloor\frac{n-(q^D-1)}{q^D}\right\rfloor.
\]

Every monomial of $\prod I_{V,i}/Y$ has degree at most $mq^{D-1}-1$.
The trace formula \eqref{eq:trace} therefore gives
\[
 \deg_XG_{V,S}\leqslant
 \left\lfloor\frac{mq^{D-1}-q^D}{q^D}\right\rfloor
 =\lfloor m/q\rfloor-1,
\]
with $G_{V,S}=0$ when $m<q$.

\emph{(ii).} Replacing $w$ by $\zeta w$, for $\zeta\in\Fq^\times$, and reindexing
the finite sum shows
\[
 G_{V,S}(\zeta X)=\zeta^{m-1}G_{V,S}(X).
\]
Thus any exponent $d$ occurring satisfies $d\equiv m-1\pmod{q-1}$.
The degree bound gives $m-q(d+1)\geqslant0$, and this congruence
makes it a multiple of $q-1$. This proves (ii).
\end{proof}

\subsection{Lattice exhaustion and interpolation limits}
\label{sec:limits}

\begin{notation}\label{not:lattice-exhaustion}
For $R>0$, put
$V_R:=\{\lambda\in\Lambda:|\lambda|\leqslant R\}$, and restrict to
$R$ for which $V_R\ne\{0\}$. Write $E_{V_R}$ and
$I_{V_R,i}$ for the polynomials in \eqref{eq:EV} and \eqref{eq:IV}, and
$b_{V_R,i,j}$ for the coefficients of $I_{V_R,i}$, as in \eqref{eq:IV}.
For $\rho>0$, the series $B(w)=\sum_jb_jw^j$ with $b_j\in F$ and
$\norm{b_j}_F\rho^j\to0$ form the complete closed-disk algebra
of radius $\rho$, with Gauss norm
\[
 \norm{B}_{F,\rho}:=\sup_j\norm{b_j}_F\rho^j.
\]
\end{notation}

\begin{lemma}\label{lem:interpolation-convergence}
\textup{(i)} When $R\to\infty$, the functions $E_{V_R}$ converge to $e$
in $\norm{\cdot}_{F,\rho}$ for every $\rho>0$, hence uniformly on
bounded disks.
\par\medskip\noindent\textup{(ii)} The polynomials $I_{V_R,i}$ converge to $h_i$ in
$\norm{\cdot}_{F,\rho}$ for every $\rho>0$, hence uniformly on
bounded disks. Moreover, $b_{V_R,i,0}\to c_i$.

\end{lemma}

\begin{proof}

\emph{(i).} Fix $\rho>0$ and take $R>\rho$. The tail of the lattice
product,
\[
 U_R(w)=\prod_{\substack{\lambda\in\Lambda\\|\lambda|>R}}
                    (1-w/\lambda),
\]
converges in the closed-disk algebra and satisfies
\[
 e=E_{V_R}U_R,\qquad
 \norm{U_R^{-1}}_{F,\rho}=1,\qquad
 \norm{U_R^{-1}-1}_{F,\rho}\leqslant\rho/R.
\]
Indeed, each factor differs from one by a term of norm at most
$\rho/R<1$, and these norms tend to zero. Thus
$\norm{E_{V_R}-e}_{F,\rho}\leqslant
(\rho/R)\norm{e}_{F,\rho}$, proving (i). The same identities give
$\norm{E_{V_R}}_{F,\rho}=\norm{e}_{F,\rho}$.

\emph{(ii).} Expanding each $(w-\lambda)^{-1}$ as a geometric series
and using Lemma~\ref{lem:period-bound}(ii) gives
\begin{equation}\label{eq:tail}
 T_{i,R}(w):=\sum_{\substack{\lambda\in\Lambda\\|\lambda|>R}}
     \frac{h_i(\lambda)}{w-\lambda},\qquad
 \norm{T_{i,R}}_{F,\rho}\leqslant M_i/R.
\end{equation}
Splitting the lattice sum in \eqref{eq:inputs} gives
$h_i=U_R I_{V_R,i}+eT_{i,R}$. This identity holds first away from
$V_R$ and then on the whole disk, since both sides are analytic.
Consequently
\[
 I_{V_R,i}-h_i=(U_R^{-1}-1)h_i-E_{V_R}T_{i,R},
\]
and therefore
\[
 \norm{I_{V_R,i}-h_i}_{F,\rho}
 \leqslant\max\left\{
       \frac\rho R\norm{h_i}_{F,\rho},
       \frac{M_i}R\norm{e}_{F,\rho}\right\}
 \longrightarrow0.
\]
Taking the coefficient of $w$ in the convergent series gives
$b_{V_R,i,0}\to c_i$.
\end{proof}

\begin{proposition}\label{prop:lattice-limit}
For $w\notin\Lambda$, put
\[
 R_S(w):=\sum_{\lambda\in\Lambda}
       \frac{\prod_{i\in S}h_i(\lambda-w)}{w-\lambda}.
\]
\begin{enumerate}[label=(\roman*)]
\item The series defining $R_S$ converges normally on closed
disks away from $\Lambda$, and
\begin{equation}\label{eq:sum-limit}
 G_{V_R,S}(E_{V_R}(w))\longrightarrow R_S(w).
\end{equation}
\item There is a unique polynomial
$G_S\in F[X]$ such that
$R_S(w)=G_S(e(w))$. The identity extends over $\Lambda$, and the apparent
singularities are removable.
\end{enumerate}
\end{proposition}

\begin{proof}

\emph{(i).} Fix a closed disk $D$ disjoint from $\Lambda$, and choose
$\rho$ with $|w|\leqslant\rho$ on $D$.
Lemma~\ref{lem:interpolation-convergence} bounds $I_{V_R,i}(w)$ uniformly
on $D$ for all sufficiently large $R$. For $\lambda\in V_R$, additivity gives
\[
 I_{V_R,i}(\lambda-w)=h_i(\lambda)-I_{V_R,i}(w).
\]
Together with Lemma~\ref{lem:period-bound}, this identity and
interpolation convergence give a constant $M_D>0$, independent of $R$
and $\lambda$, such that the finite and limiting numerators satisfy,
respectively,
\[
\begin{gathered}
 \norm{\prod_{i\in S}I_{V_R,i}(\lambda-w)}_F\leqslant M_D
 \qquad(\lambda\in V_R),\\
 \norm{\prod_{i\in S}\bigl(h_i(\lambda)-h_i(w)\bigr)}_F\leqslant M_D
 \qquad(\lambda\in\Lambda),
\end{gathered}
\]
for all sufficiently large $R$. If $|\lambda|>\rho$, then
$|w-\lambda|=|\lambda|$, and hence the corresponding summands satisfy
\[
\begin{gathered}
 \norm{\frac{\prod_{i\in S}I_{V_R,i}(\lambda-w)}{w-\lambda}}_F
 \leqslant\frac{M_D}{|\lambda|},\\
 \norm{\frac{\prod_{i\in S}\bigl(h_i(\lambda)-h_i(w)\bigr)}
              {w-\lambda}}_F
 \leqslant\frac{M_D}{|\lambda|}.
\end{gathered}
\]
For $L>\rho$, either tail with $|\lambda|>L$ thus has norm at most
$M_D/L$, uniformly in $R$. On the finite set $|\lambda|\leqslant L$,
the summands converge uniformly as $R\to\infty$ by interpolation
convergence. Letting $R\to\infty$ and then $L\to\infty$ proves
\eqref{eq:sum-limit}; the tail estimate also proves normal convergence
of $R_S$ on $D$.

\emph{(ii).} For fixed $m\geqslant q$, put
$r:=\lfloor m/q\rfloor$ and choose $w_1,\ldots,w_r$ in a sufficiently
small punctured disk such that the $e(w_j)$ are distinct; such points
exist because $e'(0)=1$. The Vandermonde matrices satisfy
\[
 \mathcal V_R:=\bigl(E_{V_R}(w_j)^d\bigr)_{j,d}
 \longrightarrow
 \mathcal V:=\bigl(e(w_j)^d\bigr)_{j,d},
\]
and
\[
 \det\mathcal V
 =\prod_{1\leqslant j<l\leqslant r}\bigl(e(w_l)-e(w_j)\bigr)\ne0,
\]
where $1\leqslant j,l\leqslant r$ and $0\leqslant d\leqslant r-1$.
Equation~\eqref{eq:sum-limit} therefore gives coefficientwise convergence
of $G_{V_R,S}$ to a polynomial $G_S\in F[X]$. For $m<q$, put $G_S:=0$,
since every $G_{V_R,S}$ is zero. For each $w\notin\Lambda$, coefficient
convergence and $E_{V_R}(w)\to e(w)$ now give
$G_{V_R,S}(E_{V_R}(w))\to G_S(e(w))$. Comparing this limit with
\eqref{eq:sum-limit} proves $R_S(w)=G_S(e(w))$.

Near a lattice point, the one potentially singular term has
numerator vanishing to order at least $m$, and the remaining series
converges normally. Thus the identity extends over $\Lambda$. Uniqueness
follows because $e$ has infinite image.
\end{proof}

\begin{notation}\label{not:lattice-coefficients}
For the limit polynomial from
Proposition~\ref{prop:lattice-limit}(ii), write
\[
 G_S(X):=\sum_{d\geqslant0}g_{S,d}X^d.
\]
Set $g_{\varnothing,-1}:=1$, and set all other coefficients
outside their polynomial ranges equal to zero.
\end{notation}

\begin{corollary}\label{cor:lattice-bounds}
One has $\deg G_S\leqslant\lfloor m/q\rfloor-1$, with a negative
bound meaning $G_S=0$. If $g_{S,d}\ne0$, then
$m=q(d+1)+(q-1)k$ for an integer $k\geqslant0$.
\end{corollary}

\begin{proof}
Proposition~\ref{prop:lattice-limit}(ii) gives coefficientwise
convergence $G_{V_R,S}\to G_S$. The coefficients excluded by
Lemma~\ref{lem:finite} vanish for every $R$ and therefore also in the
limit.
\end{proof}

\begin{remark}
The finite tensor formula will show that the support condition
in Corollary~\ref{cor:lattice-bounds} is also sufficient for nonvanishing
(Lemma~\ref{lem:support}). Proposition~\ref{prop:trichotomy} then gives
the exact degree.
\end{remark}

\subsection{Scalar and differential generating identities}

\begin{notation}\label{not:scalar-to-differential}
Write $F\nu_S$ for the one-dimensional $F$-space spanned by
$\nu_S$.
For the scalar polynomial $G_S\in F[X]$ constructed in
Proposition~\ref{prop:lattice-limit}(ii), set
\[
 \mathfrak{G}_S:=G_S\nu_S\in(F\nu_S)[X],
\]
and define $\mathfrak{g}_{S,d}:=g_{S,d}\nu_S$, so that
$\mathfrak{g}_{S,d}$ is the coefficient of $X^d$ in $\mathfrak{G}_S$.
\end{notation}

\begin{proposition}\label{prop:scalar-generating}
With $G_\varnothing(X):=X^{-1}$,
\begin{equation}\label{eq:scalar-generating}
 \sum_{0\ne\lambda\in\Lambda}
       \frac{\prod_{i=1}^s h_i(\lambda)}{w-\lambda}
 =\sum_{S\subseteq\Sigma_s}
       G_S(e(w))\prod_{i\notin S}h_i(w).
\end{equation}
\end{proposition}

\begin{proof}
The scalar partial fraction identity
\[
 \sum_{\lambda\in\Lambda}\frac1{w-\lambda}=\frac1{e(w)}
\]
follows by taking limits of
$E_{V_R}'/E_{V_R}=\sum_{\lambda\in V_R}(w-\lambda)^{-1}$ away from
$\Lambda$: here $E_{V_R}'=1$ and $E_{V_R}\to e$. For each $\lambda\in\Lambda$, write
$h_i(\lambda)=h_i(\lambda-w)+h_i(w)$. The product expansion is
\[
\begin{aligned}
\frac{\prod_{i\in\Sigma_s}h_i(\lambda)}{w-\lambda}
&=\frac{\prod_{i\in\Sigma_s}\bigl(h_i(\lambda-w)+h_i(w)\bigr)}
        {w-\lambda}\\
&=\sum_{S\subseteq\Sigma_s}
  \frac{\prod_{i\in S}h_i(\lambda-w)}{w-\lambda}
  \prod_{i\in\Sigma_s\setminus S}h_i(w).
\end{aligned}
\]
Proposition~\ref{prop:lattice-limit} identifies
the sum for each nonempty $S$, while the partial fraction identity above
handles $S=\varnothing$. The $\lambda=0$ term vanishes because
$h_i(0)=0$ for every $i$, so summation over all lattice points agrees
with the original sum over nonzero lattice points. This proves
\eqref{eq:scalar-generating}.
\end{proof}

\begin{definition}\label{def:mc-pellarin-perkins}
The multicharacter Pellarin--Perkins generating function is
\[
 \mathscr C_{\varphi}^s(w):=\sum_{0\ne\lambda\in\Lambda}
       \frac{\Hcal_1(\lambda)\boxtimes\cdots\boxtimes\Hcal_s(\lambda)}
            {w-\lambda}.
\]
\end{definition}

\begin{samepage}
\begin{theorem}
\label{thm:main}
\textup{(i)} For every nonempty $S$ there is a unique polynomial
$\mathfrak{G}_S(X)\in(F\nu_S)[X]$ characterized by
\begin{equation}\label{eq:canonical}
 \mathfrak{G}_S(e(w))
 =\sum_{\lambda\in\Lambda}
       \frac{\prodext_{i\in S}\Hcal_i(\lambda-w)}{w-\lambda}.
\end{equation}

\par\medskip\noindent\textup{(ii)} With $\mathfrak{G}_\varnothing(X):=X^{-1}$, the generating identity
holds meromorphically in $w$:
\begin{equation}\label{eq:generating}
 \mathscr C_\varphi^s(w)
 =\sum_{S\subseteq\Sigma_s}
       \mathfrak{G}_S(e(w))\boxtimes\prodext_{i\notin S}\Hcal_i(w).
\end{equation}

\end{theorem}
\end{samepage}

\begin{proof}

\emph{(i).} Proposition~\ref{prop:lattice-limit} gives
$R_S(w)=G_S(e(w))$ with $G_S\in F[X]$ and removable lattice
singularities. Multiplication by $\nu_S$ gives \eqref{eq:canonical} with
$\mathfrak G_S\in(F\nu_S)[X]$. Since $e$ has infinite image,
\eqref{eq:canonical} determines $\mathfrak G_S$ uniquely.
Its right-hand side involves only the generating differentials
$\Hcal_i$, so $\mathfrak G_S$ is independent of the choice of $\nu$.

\emph{(ii).} Multiply \eqref{eq:scalar-generating} by $\nu_{\Sigma_s}$.
Using $\mathfrak{G}_S=G_S\nu_S$ and $\Hcal_i=h_i\nu_i$, the term
indexed by $S$ becomes
$\mathfrak{G}_S(e(w))\boxtimes\prodext_{i\notin S}\Hcal_i(w)$.
This proves \eqref{eq:generating}.
\end{proof}

\section{Frobenius recursion and explicit correction coefficients}\label{sec:coeff-formula}

The Frobenius--Hasse identity below gives a recursion for the scalar
correction coefficients. After normalization by the shtuka differentials,
the recursion terminates and can be indexed by depth tuples. This yields
the finite tensor formula in Theorem~\ref{thm:formula}, which proves
rationality and determines the support and degrees.

\subsection{A Frobenius--Hasse recursion for every coefficient}
\label{sec:FH}
\begin{definition}[Hasse derivative]
For a Laurent polynomial $G$, the Hasse derivative $D_X^{[a]}G$ is 
defined by the formal expansion in $F(X)[[v]]$
\[ G(X+v)=\sum_{a\geqslant0}D_X^{[a]}G(X)v^a.\]
\end{definition}

For the next lemma, suppose that the field $F$
introduced in Section~\ref{sec:preliminaries} has a continuous
automorphism $(1)$ extending $c\mapsto c^q$ on $\Cinf$, and suppose that
\begin{equation}\label{eq:eigenmaps}
 h_i(\lambda)^{(1)}=f_i h_i(\lambda)\quad(\lambda\in V),\qquad f_i\ne0.
\end{equation}
In the geometric application, coefficient Frobenius satisfies
this hypothesis and is continuous, with
$\norm{u^{(1)}}_F=\norm{u}_F^q$.
For polynomials in $w$, write
$\tau_w(I_{V,i})(w)=I_{V,i}^{(1)}(w^q)$, where the superscript on the
right twists only coefficients.

\begin{lemma}\label{lem:FH}
Assume \eqref{eq:eigenmaps} for the finite interpolation data.
Put $f_S:=\prod_{i\in S}f_i$, $b_{V,T,0}:=\prod_{i\in T}b_{V,i,0}$,
and $c_T:=\prod_{i\in T}c_i$. Here $b_{V,i,0}$ is as in
\eqref{eq:IV}, and $c_i=C_{\varphi,i}/\nu_i$ is as in
Notation~\ref{not:scalar-trivialization}.
\begin{enumerate}[label=(\roman*)]
\item With $G_{V,\varnothing}:=X^{-1}$,
\begin{equation}\label{eq:FH-finite}
 G_{V,S}^{(1)}(X^q)
 =f_S D_X^{[q-1]}
       \left(\sum_{T\subseteq S}b_{V,S\setminus T,0}
                     X^{m-|T|}G_{V,T}(X)\right).
\end{equation}
\item With $G_\varnothing(X):=X^{-1}$, the lattice-limit identity is
\begin{equation}\label{eq:FH-lattice}
 G_S^{(1)}(X^q)
 =f_S D_X^{[q-1]}
       \left(\sum_{T\subseteq S}c_{S\setminus T}
                     X^{m-|T|}G_T(X)\right).
\end{equation}
\end{enumerate}
\end{lemma}

\begin{proof}
The expression $\tau_w(I_{V,i})-f_i I_{V,i}$ vanishes on $V$
and has degree at most $q^D$, so it is a scalar multiple of $E_V$.
Its derivative at zero is $-f_i b_{V,i,0}$, and hence
\[
 \tau_w(I_{V,i})-f_i I_{V,i}
 =-f_i b_{V,i,0}E_V.
\]
Comparing coefficients of $w^{q^D}$ gives
\begin{equation}\label{eq:top-twist}
 b_{V,i,D-1}^{(1)}=-f_i b_{V,i,0}\alpha.
\end{equation}

\emph{(i).} Put
\[
 \Phi_S(X):=\sum_{T\subseteq S}b_{V,S\setminus T,0}
                     X^{m-|T|}G_{V,T}(X).
\]
Write $X:=E_V(w)$ and, for each $\lambda\in V$, put
$u_{i,\lambda}:=I_{V,i}(\lambda-w)$. Expanding the product gives
\[
 \Phi_S(X)
 =\sum_{\lambda\in V}\frac{\prod_{i\in S}(u_{i,\lambda}+b_{V,i,0}X)}{w-\lambda}.
\]
This uses $\sum_{\lambda\in V}(w-\lambda)^{-1}=X^{-1}$ for the empty
subset.

Work in $F(w)[[v]]$. Since
$E_V(w+v)\equiv X+v\pmod{v^q}$, one has
$\Phi_S(E_V(w+v))\equiv \Phi_S(X+v)\pmod{v^q}$.
Also $I_{V,i}(\lambda-w-v)\equiv u_{i,\lambda}-b_{V,i,0}v\pmod{v^q}$, so every
numerator on the right is unchanged modulo $v^q$.
Taking the coefficient of $v^{q-1}$ therefore gives
\[
 D_X^{[q-1]}\Phi_S(X)
 =\sum_{\lambda\in V}\frac{\prod_{i\in S}(u_{i,\lambda}+b_{V,i,0}X)}{(w-\lambda)^q},
\]
since $(-1)^{q-1}=1$ in characteristic $p$.

Evaluating the identity for
$\tau_w(I_{V,i})-f_i I_{V,i}$ at $\lambda-w$ gives
\[
 I_{V,i}^{(1)}((\lambda-w)^q)=f_i(u_{i,\lambda}+b_{V,i,0}X).
\]
Apply the map $R(w)\mapsto R^{(1)}(w^q)$ to the finite
sum identity in Lemma~\ref{lem:finite}. It sends each scalar
$\lambda\in V$ to $\lambda^q$ and fixes the geometric variables. Because
$E_V^{(1)}(w^q)=E_V(w)^q$, it gives
\[
 \begin{aligned}
 G_{V,S}^{(1)}(X^q)
 &=\sum_{\lambda\in V}
     \frac{\prod_{i\in S}I_{V,i}^{(1)}(\lambda^q-w^q)}{w^q-\lambda^q}\\
 &=f_S\sum_{\lambda\in V}
     \frac{\prod_{i\in S}(u_{i,\lambda}+b_{V,i,0}X)}{(w-\lambda)^q}.
 \end{aligned}
\]
Comparing the last two displays proves \eqref{eq:FH-finite} after
$X=E_V(w)$. Substitution by the nonconstant polynomial $E_V(w)$ is
injective on $F(X)$, so the identity holds in $X$.

\emph{(ii).} Apply (i) with $V=V_R$.
Lemma~\ref{lem:interpolation-convergence}(ii) gives
$b_{V_R,i,0}\to c_i$, while Proposition~\ref{prop:lattice-limit} gives
coefficientwise convergence $G_{V_R,T}\to G_T$ for every nonempty $T$.
The degrees are bounded independently of $R$ by Lemma~\ref{lem:finite}(i),
and the empty-subset term is the fixed Laurent polynomial $X^{-1}$.
On these bounded-degree expressions, coefficient Frobenius, finite
products and sums, and the fixed Hasse derivative preserve coefficientwise
limits. Taking the limit in \eqref{eq:FH-finite} proves
\eqref{eq:FH-lattice}.
\end{proof}

\begin{lemma}\label{lem:coefficient-recursion}
Recall the coefficient conventions of
Notation~\ref{not:lattice-coefficients}. For every $d\geqslant0$, put
$h:=m-q(d+1)$. Then
\begin{equation}\label{eq:coefficient-recursion}
 g_{S,d}^{(1)}
 =f_S\sum_{T\subseteq S}c_{S\setminus T}
                     g_{T,\,|T|-h-1}.
\end{equation}
\end{lemma}

\begin{proof}
Take the coefficient of $X^{qd}$ in \eqref{eq:FH-lattice}.
Before Hasse differentiation, the contributing exponent is
$qd+q-1$. The congruence
$(1+v)^{qd+q-1}\equiv(1+v)^{q-1}\pmod{v^q}$ gives
$\binom{qd+q-1}{q-1}=1$ in characteristic $p$. For a term
$X^{m-|T|}g_{T,j}X^j$, the required exponent condition is
\[
 m-|T|+j=qd+q-1,
\]
or equivalently,
\[
 j=qd+q-1-m+|T|=|T|-h-1.
\]
Summing over $T$ gives \eqref{eq:coefficient-recursion}.
\end{proof}

\subsection{The normalized coefficient recursion}

\begin{notation}\label{not:normalized-coefficients}
\leavevmode 
\begin{enumerate}
      \item Write $\Delta_{\varphi,i}$ for the normalized shtuka differential on the
$i$-th geometric factor, and put $\delta_i:=\Delta_{\varphi,i}/\nu_i$. Define
\[
 \gamma_i:=\frac{(f_i \delta_i)^{(-1)}}{\delta_i}.
\]
The identity $\delta_i^{(1)}=f_ic_i$ gives
$\gamma_i^{(1)}=f_i\delta_i/\delta_i^{(1)}=\delta_i/c_i$.

\item Whenever $m=q(d+1)+(q-1)k$ for some $d,k\geqslant0$, define
\[ 
P_{S,k}:=g_{S,d}/\prod_{i\in S}\delta_i.
\]
\end{enumerate}
\end{notation}
\begin{remark}
For fixed $S$, the relation $m=q(d+1)+(q-1)k$ determines either
index from the other. The index $k$ records the progress of the recursion:
it strictly decreases at every step and terminates at zero.
\end{remark}

\begin{proposition}\label{prop:recursion}
For every nonempty $S$ with $q\mid m$, one has $P_{S,0}=1$.
For $d\geqslant0$ and $k\geqslant1$ satisfying
$m=q(d+1)+(q-1)k$,
\begin{equation}\label{eq:normalized-recursion}
 P_{S,k}=\sum_{j=0}^{k-1}
       \ \sum_{\substack{T\subseteq S\\|T|=qk-j}}
                  \left(\prod_{i\in T}\gamma_i\right)P_{T,j}^{(-1)}.
\end{equation}
\end{proposition}

\Needspace{16\baselineskip}
\begin{proof}
For $k=0$, write $m=qr$ with $r\geqslant1$.
Let $D=\dim_{\Fq}V_R$, let $\alpha$ be the leading coefficient of
$E_{V_R}$, and let $L_{S,V_R}$ be the coefficient of $X^{r-1}$ in
$G_{V_R,S}$. In \eqref{eq:resolvent}, a term containing $X^{r-1}$
has $T$-degree at least $rq^D-1$, which is also the upper bound for the
$Y$-degree of $\prod_{i\in S}I_{V_R,i}(Y)/Y$.
The coefficient of $X^{r-1}T^{rq^D-1}$ in \eqref{eq:resolvent} is
$(-1)^{r-1}\alpha^{-r}$. The trace formula \eqref{eq:trace} therefore gives
\[
 L_{S,V_R}=(-1)^r\alpha^{-r}\prod_{i\in S}b_{V_R,i,D-1}.
\]
Twisting and substituting \eqref{eq:top-twist} yield
\[
 L_{S,V_R}^{(1)}
 =(-1)^{r+m}\alpha^{m-qr}\prod_{i\in S}f_i b_{V_R,i,0}
 =\prod_{i\in S}f_i b_{V_R,i,0}.
\]
Here $m=qr$, and $(-1)^{r+m}=1$ in characteristic $p$. Since
$b_{V_R,i,0}\to c_i$, continuity of twisting yields
\[
 g_{S,r-1}^{(1)}=\prod_{i\in S}f_ic_i,
\]
and hence $P_{S,0}=1$.
For $k>0$, put $h:=(q-1)k$. The $T=\varnothing$ term in
\eqref{eq:coefficient-recursion} is zero because
$g_{\varnothing,-h-1}=0$. A nonzero term indexed by $T\ne\varnothing$
requires, by Corollary~\ref{cor:lattice-bounds}, an integer
$j\geqslant0$ such that
\[
 |T|=q(|T|-h)+(q-1)j.
\]
Thus $|T|=qk-j$, and the coefficient index $|T|-h-1$ equals
$k-j-1$. This index is nonnegative if and only if $j<k$.

Since $\delta_i^{(1)}=f_ic_i$, dividing
\eqref{eq:coefficient-recursion} by
$\prod_{i\in S}\delta_i^{(1)}=f_S\prod_{i\in S}c_i$ gives
\[
 P_{S,k}^{(1)}
 =\sum_{j=0}^{k-1}\ \sum_{\substack{T\subseteq S\\|T|=qk-j}}
       \left(\prod_{i\in T}\frac{\delta_i}{c_i}\right)P_{T,j}.
\]
Replacing $\delta_i/c_i$ by $\gamma_i^{(1)}$ and taking one inverse twist
gives \eqref{eq:normalized-recursion}. Since $j<k$, the recurrence
determines the $P_{S,k}$ by induction on $k$.
\end{proof}

\subsection{Depth tuples and expansion chains}
\label{sec:formula}

\begin{notation}\label{not:depth}
Let $S$ be a nonempty finite set. A \emph{depth tuple} on $S$ is
a tuple $(n_i)_{i\in S}$ of nonnegative integers. Its \emph{weight} is
$\sum_{i\in S}q^{-n_i}$. For an integer $d\geqslant0$, let
$\mathbb E_{S,d}$ be the set of depth tuples of weight $q(d+1)$.

\end{notation}
The \emph{depth profile} of a tuple $(n_i)$ is the sequence
$x_j:=\#\{i\in S:n_i=j\}$ for $j\geqslant0$.
\begin{lemma}\label{lem:depth-profile}
For $(n_i)\in\mathbb E_{S,d}$, its depth profile satisfies
\[
 \sum_{j\geqslant0}x_j=|S|,\qquad
 \sum_{j\geqslant0}q^{-j}x_j=\sum_{i\in S}q^{-n_i}=q(d+1).
\] 
\end{lemma}
\begin{lemma}\label{lem:depth-existence}
For every nonempty $S$ and every $d\geqslant0$, put $m:=|S|$.
The set $\mathbb{E}_{S,d}$ is nonempty if and only if
$m=q(d+1)+(q-1)k$ for some integer $k\geqslant0$.
\end{lemma}

\begin{proof}
If a depth tuple satisfies
$\sum_{i\in S}q^{-n_i}=q(d+1)$, then $q(d+1)\leqslant m$. Multiplying by
$q^{\max n_i}$ and reducing modulo $q-1$ gives
$m\equiv q(d+1)\pmod{q-1}$, so $m=q(d+1)+(q-1)k$ for some $k\geqslant0$.

Conversely, suppose that $m=q(d+1)+(q-1)k$ with $k\geqslant0$.
For $k=0$, take $n_i=0$ for every $i\in S$; this has $m=q(d+1)$ indices
and total weight $q(d+1)$. For $k\geqslant1$, take $q(d+1)-1$ indices of
depth zero, $q-1$ indices of depth $j$ for each $1\leqslant j\leqslant k-1$,
and $q$ indices of depth $k$. The number of indices is
\[
 q(d+1)-1+(q-1)(k-1)+q=q(d+1)+(q-1)k=m,
\]
and the total weight is
\[
 q(d+1)-1+(q-1)\sum_{j=1}^{k-1}q^{-j}+q^{1-k}=q(d+1).
\]
Thus the required depth tuple exists.
\end{proof}

\begin{definition}\label{def:expansion-chains}
Fix a nonempty finite set $S$, put $m:=|S|$, and let
$d,k\geqslant0$ be integers satisfying $m=q(d+1)+(q-1)k$.
Put $T_0:=S$ and $k_0:=k$. An \emph{expansion chain}
consists of an integer $t\geqslant0$ and, when $t>0$, subsets
\[
 T_0=S\supseteq T_1\supseteq\cdots\supseteq T_t,
\]
together with integers 
\[ 
 k = k_0 > k_1>\ldots > k_t =0 \]
  satisfying  $
 |T_j|=qk_{j-1}-k_j\quad(1\leqslant j\leqslant t)$.
Let $\mathfrak C_{S,d,k}$ denote the set of all expansion chains.
When $t=0$, we require $k=0$, and the chain consists only of
the initial pair $(T_0,k_0)=(S,0)$.
For formulas involving cardinalities beyond the chain length, put
$T_j:=\varnothing$ for $j>t$.
\end{definition}
\begin{remark}\label{rem:expansion-chains}
      An expansion chain records successive choices in
\eqref{eq:normalized-recursion}: at step $j$, the summand indexed by
$(T_j,k_j)$ replaces $P_{T_{j-1},k_{j-1}}$. The terminal coefficient is
$P_{T_t,0}=1$; when $t=0$, it is the initial coefficient $P_{S,0}$.
\end{remark}

For
$\mathfrak c\in\mathfrak C_{S,d,k}$, write $t(\mathfrak c)$ for its length
and define the $i$-th depth of $\mathfrak c$ by
\[
 \operatorname{dep}_i(\mathfrak c):=\#\{j:1\leqslant j\leqslant t(\mathfrak c)
 \text{ and }i\in T_j\}\qquad(i\in S).
\]
The induced depth tuples define a map
\[
 \operatorname{dep}:\mathfrak C_{S,d,k}\longrightarrow
 \mathbb Z_{\geqslant0}^S,\qquad
 \mathfrak c\longmapsto\bigl(\operatorname{dep}_i(\mathfrak c)\bigr)_{i\in S}.
\] 

\begin{lemma}\label{lem:chain-combinatorics}
For $m=q(d+1)+(q-1)k$ with $d,k\geqslant0$, the following hold.
\begin{enumerate}[label=(\roman*)]
\item The map $\operatorname{dep}$ takes values in $\mathbb E_{S,d}$; that is,
$\operatorname{dep}(\mathfrak C_{S,d,k})\subseteq\mathbb E_{S,d}$.
Moreover, for every $\mathfrak c\in\mathfrak C_{S,d,k}$, its coordinates
satisfy $\operatorname{dep}_i(\mathfrak c)\leqslant t(\mathfrak c)\leqslant k$ for every
$i\in S$.
\item The resulting map
\[
 \operatorname{dep}:\mathfrak C_{S,d,k}\longrightarrow\mathbb E_{S,d}
\]
is a bijection.
\item For $\mathfrak c\in\mathfrak C_{S,d,k}$, the depth profile
$(x_j)_{j\geqslant0}$ of $\operatorname{dep}(\mathfrak c)$ satisfies
\begin{equation}\label{eq:depth-profile}
 |T_j|=\sum_{a\geqslant j}x_a,\qquad
 x_j=|T_j|-|T_{j+1}|\qquad(j\geqslant0).
\end{equation}
In particular, the depth profile determines every cardinality $|T_j|$,
and the sequence $(|T_j|)_{j\geqslant1}$, together with the fixed integer
$|S|$, determines the depth profile.
\end{enumerate}
\end{lemma}

\begin{proof}
\emph{(i).} Fix $\mathfrak c\in\mathfrak C_{S,d,k}$, and put
$t:=t(\mathfrak c)$ and $(n_i):=\operatorname{dep}(\mathfrak c)$. The chain
subsets are nonempty because
\[
 |T_j|=qk_{j-1}-k_j>(q-1)k_{j-1}\geqslant0
 \qquad(1\leqslant j\leqslant t).
\]
The cardinality relation and the definition of $n_i$ give
\[
 \sum_{j=1}^t|T_j|q^{-j}
 =\sum_{j=1}^t\left(k_{j-1}q^{1-j}-k_jq^{-j}\right)
 =k_0-k_tq^{-t}=k
\]
and
\[
 \sum_{j=1}^t|T_j|q^{-j}
 =\sum_{i\in S}\sum_{j=1}^{n_i}q^{-j}
 =\frac1{q-1}\sum_{i\in S}\left(1-q^{-n_i}\right).
\]
Therefore
\begin{equation}\label{eq:depth-weight}
 k=\sum_{j=1}^t|T_j|q^{-j}
 =\frac1{q-1}\sum_{i\in S}\left(1-q^{-n_i}\right),
\end{equation}
and
\begin{equation}\label{eq:depth-weight2} 
 \sum_{i\in S}q^{-n_i}=m-(q-1)k=q(d+1).
\end{equation}
The weight identity places $(n_i)$ in $\mathbb E_{S,d}$.
Each depth is at most the chain length $t$, and the strictly decreasing
integers $k=k_0>\cdots>k_t=0$ give $t\leqslant k$.

\emph{(ii).} For $(n_i)\in\mathbb E_{S,d}$, put
$\ell:=\max_i n_i$ and define
\[
 T_j:=\{i\in S:n_i\geqslant j\},\qquad u_j:=|T_j|
       \qquad(1\leqslant j\leqslant\ell).
\]
The sets are nested and nonempty. The weight equation gives
\[
 k=\frac1{q-1}\sum_{i\in S}(1-q^{-n_i})
   =\sum_{a=1}^{\ell}u_aq^{-a}.
\]
Set $k_0:=k$ and $k_j:=qk_{j-1}-u_j$. These are integers, and induction
gives
\[
 k_j=\sum_{a=j+1}^{\ell}u_aq^{j-a}\geqslant0,
 \qquad k_\ell=0.
\]
Since $u_a\leqslant u_j$ for $a\geqslant j$, for each
$1\leqslant j\leqslant\ell$ we have
\[
 0<k_{j-1}=\sum_{a=j}^{\ell}u_aq^{j-1-a}
 <\sum_{a=j}^{\infty}u_jq^{j-1-a}=\frac{u_j}{q-1}.
\]
Thus $k_j<k_{j-1}$. If $\ell=0$, the weight equation gives $k=0$,
and the construction yields the zero-length chain.

For $\ell>0$, the cardinalities satisfy $|T_j|=qk_{j-1}-k_j$.
Moreover,
\[
 |T_j|=q\bigl((k_{j-1}-k_j-1)+1\bigr)+(q-1)k_j,
\]
where $k_{j-1}-k_j-1\geqslant0$. Hence every pair $(T_j,k_j)$ has
the degree index required by Notation~\ref{not:normalized-coefficients},
and the chain occurs in the recursion.

The reconstructed chain has depth tuple $(n_i)$. Conversely, nested
nonempty chain subsets are recovered from their depths by
$T_j=\{i:n_i\geqslant j\}$, their number is $\max_i n_i$, and the
cardinality relations recover all $k_j$. The two constructions are
inverse, proving the bijection.

\emph{(iii).} Write $(n_i)=\operatorname{dep}(\mathfrak c)$.
The nested subsets satisfy $T_j=\{i\in S:n_i\geqslant j\}$, with
$T_0=S$ and $T_j=\varnothing$ for $j>t$. Counting by depth gives
$|T_j|=\sum_{a\geqslant j}x_a$, and subtracting consecutive counts gives
$x_j=|T_j|-|T_{j+1}|$. These are \eqref{eq:depth-profile}; they also give
the two stated recovery formulas.

\end{proof}

\begin{example}\label{ex:k4-profiles}
Take $q=2$, $d=0$, and $k=4$, so $|S|=6$.
By Lemma~\ref{lem:chain-combinatorics}(i)--(ii),
every tuple in $\mathbb E_{S,0}$ satisfies $n_i\leqslant k=4$, so
$x_j=0$ for $j>4$. The depth profile therefore satisfies
\[
\sum_{j=0}^4 x_j=6,\qquad
x_0+\frac{x_1}{2}+\frac{x_2}{4}+\frac{x_3}{8}+\frac{x_4}{16}=2.
\]
Equivalently,
\[
15x_0+7x_1+3x_2+x_3=26,\qquad
x_4=6-x_0-x_1-x_2-x_3.
\]
The nonnegative integer solutions give exactly the following five depth profiles.
\begin{center}
\begin{tabular}{c|c|c|c}
\hline
No. & Depth profile \((x_0,x_1,x_2,x_3,x_4)\) & \(k\)-sequence & Chain count\\
\hline
1 & \((0,2,4,0,0)\) & \(4>2>0\) & \(15\)\\[1mm]
2 & \((0,3,1,2,0)\) & \(4>2>1>0\) & \(60\)\\[1mm]
3 & \((1,0,3,2,0)\) & \(4>3>1>0\) & \(60\)\\[1mm]
4 & \((1,1,0,4,0)\) & \(4>3>2>0\) & \(30\)\\[1mm]
5 & \((1,1,1,1,2)\) & \(4>3>2>1>0\) & \(360\)\\
\hline
\end{tabular}
\end{center}
The last column counts the depth tuples with each profile:
$6!/(x_0!x_1!x_2!x_3!x_4!)$. By the bijection in
Lemma~\ref{lem:chain-combinatorics}(ii), these counts give
$15+60+60+30+360=525$ expansion chains in total.
\end{example}

\subsection{The finite tensor formula}

\leavevmode
Write $\Gamma_i^{[n]}$ for the pullback of the function
$\Gamma^{[n]}$ from Notation~\ref{not:Gamma} to the $i$-th factor.
Iterating \eqref{eq:normalized-recursion} assigns a product of these
functions to each expansion chain.
\begin{lemma}\label{lem:chain-contribution}
Fix integers $d,k\geqslant0$ satisfying
$m=q(d+1)+(q-1)k$. Expand $P_{S,k}$ by iterating
\eqref{eq:normalized-recursion}, starting with $(T_0,k_0)=(S,k)$ and
continuing until every terminal recursion index is zero. The following
assertions hold.
\begin{enumerate}[label=(\roman*)]
\item For every $\mathfrak c\in\mathfrak C_{S,d,k}$, put
$t:=t(\mathfrak c)$ and $(n_i):=\operatorname{dep}(\mathfrak c)$. The
complete summand indexed by $\mathfrak c$, including its terminal factor,
is
\begin{equation}\label{eq:chain-contribution}
\begin{aligned}
 \left(\prod_{j=1}^t\prod_{i\in T_j}\gamma_i^{(-(j-1))}\right)
 P_{T_t,0}^{(-t)}
 &=\prod_{j=1}^t\prod_{i\in T_j}\gamma_i^{(-(j-1))}\\
 &=\prod_{i\in S}\Gamma_i^{[n_i]},
\end{aligned}
\end{equation}
with the empty product equal to $1$ when $t=0$.
\item The normalized coefficient has the expansion
\[
\begin{aligned}
 P_{S,k}
 &=\sum_{\mathfrak c\in\mathfrak C_{S,d,k}}
       \prod_{i\in S}\Gamma_i^{[\operatorname{dep}_i(\mathfrak c)]}\\
 &=\sum_{(n_i)\in\mathbb E_{S,d}}
       \prod_{i\in S}\Gamma_i^{[n_i]}.
\end{aligned}
\]
Moreover, $\mathbb E_{S,d}$ is finite, and every
$(n_i)\in\mathbb E_{S,d}$ satisfies $n_i\leqslant k$ for every $i\in S$.
\end{enumerate}
\end{lemma}

\begin{proof}
\emph{(i).} Fix a chain of length $t$. The factor
$\prod_{i\in T_j}\gamma_i$ introduced at step $j$ undergoes $j-1$
inverse twists in the full expansion, and the terminal factor undergoes
$t$ inverse twists. These twists give the first expression in
\eqref{eq:chain-contribution}. Its terminal coefficient is one by
Proposition~\ref{prop:recursion}, since $|T_t|=qk_{t-1}$ for $t>0$
and $|T_0|=q(d+1)$ for $t=0$.

Because the subsets are nested, $i\in T_j$ exactly when
$1\leqslant j\leqslant n_i$. Thus the remaining product is
\[
 \prod_{j=1}^t\prod_{i\in T_j}\gamma_i^{(-(j-1))}
 =\prod_{i\in S}\prod_{r=0}^{n_i-1}\gamma_i^{(-r)}
 =\prod_{i\in S}\Gamma_i^{[n_i]}.
\]

\emph{(ii).} We induct on $k$. At $k=0$, the only chain has
length zero and contributes $P_{S,0}=1$. For $k>0$, a first choice
$(T,j)$ in \eqref{eq:normalized-recursion} satisfies $j<k$ and
\[
 |T|=qk-j=q\bigl((k-j-1)+1\bigr)+(q-1)j.
\]
The induction hypothesis therefore expands $P_{T,j}$ over
$\mathfrak C_{T,k-j-1,j}$. Prepending $(T,j)$ identifies these chains
with the chains in $\mathfrak C_{S,d,k}$ having that first step.
The inverse twist and the factor $\prod_{i\in T}\gamma_i$ in the
recursion give exactly the contribution in part~(i). Summing over
the first choices proves the chain expansion.

Lemma~\ref{lem:chain-combinatorics} replaces the chains bijectively by
depth tuples and bounds every depth by $k$. Hence
$\mathbb E_{S,d}\subseteq\{0,\ldots,k\}^S$ is finite, and the second
expansion follows.
\end{proof}

\begin{notation}\label{not:differential-spaces}
Use the projections $p_i:{\mathcal X}_{\Cinf}^S\to {\mathcal X}_{\Cinf}$ to define
\[
 \Dcal_S:=
 H^0\left(\operatorname{Spec}\Cinf({\mathcal X}^S),
          \bigotimes_{i\in S}p_i^*\Omega^1_{{\mathcal X}_{\Cinf}/\Cinf}\right).
\]
We regard this space as a space of rational external
differentials on ${\mathcal X}_{\Cinf}^{\Sigma_s}$ by pullback. All geometric differentials
after base change are relative to $\Cinf$.

Multiplication by $\nu_S$ identifies $\Cinf({\mathcal X}^S)$ with $\Dcal_S$.
After extending scalars along $\Cinf({\mathcal X}^S)\subset F$,
this gives $\Dcal_S\subset F\nu_S$.
\end{notation}

\begin{lemma}\label{lem:external-tensor-independence}
For every nonempty $S$, put
$M:=\Gamma(U_{\Cinf},\Omega^1(-V_\varphi))$ and
$\mathscr M_S:=\bigotimes_{i\in S}M$, with the tensor product over
$\Cinf$. The external products $\prodext_{i\in S}\mathfrak{B}_i^{[n_i]}$, indexed by
depth tuples on $S$, form a $\Cinf$-basis of $\mathscr M_S$ and remain
$\Cinf$-linearly independent in $\Dcal_S$ under localization at the generic
point.
\end{lemma}

\begin{proof}
Put $\mathscr A_S:=\Gamma(U_{\Cinf}^S,\mathcal O)$.
The affine external-product
identification makes $\mathscr M_S$ a
projective $\mathscr A_S$-module of rank one. The ring $\mathscr A_S$ is a domain because $U_{\Cinf}^S$ is
integral. The projective module $\mathscr M_S$ is therefore torsion-free,
so localization at the generic point is injective:
\[
 \mathscr M_S\hookrightarrow
 \operatorname{Frac}(\mathscr A_S)\otimes_{\mathscr A_S}\mathscr M_S
 \simeq\Dcal_S.
\]
The last identification is induced by
$\Omega^1(-V_\varphi)\hookrightarrow\Omega^1$, which is an
isomorphism at the generic point.
Proposition~\ref{prop:basis} supplies a $\Cinf$-basis of $M$, so
the products indexed by depth tuples form a $\Cinf$-basis of
$\mathscr M_S$ and remain linearly independent in $\Dcal_S$.
\end{proof}

\begin{samepage}
\begin{theorem}\label{thm:formula}
\textup{(i)} For every integer $d\geqslant0$, the coefficient of $X^d$ in
$\mathfrak{G}_S$ is given by
\begin{equation}\label{eq:coefficient-formula}
\mathfrak{g}_{S,d}
 =\sum_{ (n_i)\in \mathbb{E}_{S,d} }
                  \prodext_{i\in S}\mathfrak{B}_{i}^{[n_i]}.
\end{equation}
The sum is finite. Its terms are distinct elements of an external tensor
basis, and each has coefficient one.

\par\medskip\noindent\textup{(ii)} For every $d\geqslant0$, one has
\[
 \mathfrak g_{S,d}\in
 \Gamma\!\left(U_{\Cinf}^S,
       \bigotimes_{i\in S}p_i^*
          \bigl(\Omega^1_{U_{\Cinf}/\Cinf}(-V_\varphi|_{U_{\Cinf}})\bigr)
       \right).
\]
Thus it is a regular external differential and vanishes along the
pullback of $V_\varphi|_{U_{\Cinf}}$ on each factor, with the
multiplicities prescribed by that divisor.

\end{theorem}
\end{samepage}
\begin{proof}
\emph{(i).} If there is no integer $k\geqslant0$ with
$m=q(d+1)+(q-1)k$, Lemma~\ref{lem:depth-existence} makes
$\mathbb E_{S,d}$ empty and Corollary~\ref{cor:lattice-bounds} gives
$\mathfrak g_{S,d}=0$. Thus \eqref{eq:coefficient-formula} holds.

Now suppose $m=q(d+1)+(q-1)k$ with $k\geqslant0$.
By the definition of $P_{S,k}$ and
Lemma~\ref{lem:chain-contribution}(ii), together with
$\Delta_{\varphi,i}\Gamma_i^{[n_i]}=\mathfrak{B}_i^{[n_i]}$ from the proof of
Proposition~\ref{prop:basis}, we obtain
\[
 \mathfrak g_{S,d}
 =\left(\prod_{i\in S}\delta_i\right)P_{S,k}\nu_S
 =\sum_{(n_i)\in\mathbb E_{S,d}}
       \prodext_{i\in S}
       \left(\Delta_{\varphi,i}\Gamma_i^{[n_i]}\right)
 =\sum_{(n_i)\in\mathbb E_{S,d}}
       \prodext_{i\in S}\mathfrak{B}_i^{[n_i]}.
\]
The same lemma gives finiteness and the depth bound $n_i\leqslant k$.
By Lemma~\ref{lem:external-tensor-independence}, these products
belong to an external tensor basis and remain linearly independent in
$\Dcal_S$. The bijection in Lemma~\ref{lem:chain-combinatorics} ensures
that each occurs once, proving the remaining assertions in part~(i).

\emph{(ii).} Each $\mathfrak{B}^{[n]}$ is a section of
$\Omega^1_{U_{\Cinf}/\Cinf}(-V_\varphi|_{U_{\Cinf}})$ by
Proposition~\ref{prop:basis}. Every term of the finite sum
\eqref{eq:coefficient-formula}, and hence the coefficient itself,
belongs to the external tensor space in part (ii).
\end{proof}

\begin{corollary}\label{cor:rational-coefficients}
 For every nonempty $S$, we have
\[
 G_S\in\Cinf({\mathcal X}^S)[X],
 \qquad
 \mathfrak G_S\in\Dcal_S[X].
\]
\end{corollary}

\begin{proof}
Theorem~\ref{thm:formula}(i) expresses every coefficient of
$\mathfrak G_S$ as a finite sum of rational external differentials,
so $\mathfrak G_S\in\Dcal_S[X]$. Dividing by $\nu_S$ gives
$G_S\in\Cinf({\mathcal X}^S)[X]$ by
Notation~\ref{not:differential-spaces}.

\end{proof}
\begin{remark}
      This rationality completes the construction of the canonical corrections
in Definition~\ref{def:correction-polynomials}. Together with Theorem~\ref{thm:main}, it proves
Theorem~\ref{thm:intromain}(i)--(ii). The constancy criterion in
part~(iii) follows from Proposition~\ref{prop:trichotomy} and
Corollary~\ref{cor:simultaneous-constancy} below.
\end{remark}

\subsection{Classification, degrees, and leading forms}\label{sec:support-degrees}
\begin{lemma}\label{lem:support}
 For every $d\geqslant0$, ${\mathfrak{g}_{S,d}}\ne 0$ if and only if 
\begin{equation}\label{eq:support}
 m=q(d+1)+(q-1)k\quad\mbox{for some integer }k\geqslant0.
\end{equation}
\end{lemma}
\begin{proof}
By Theorem~\ref{thm:formula}(i), the coefficient is a sum of distinct
external tensor basis elements, each with coefficient one. It is
therefore nonzero if and only if $\mathbb E_{S,d}$ is nonempty.
Lemma~\ref{lem:depth-existence} gives precisely \eqref{eq:support}.
\end{proof}
Recall that $\Nq=q\mathbb Z_{\geqslant2}
 +(q-1)\mathbb Z_{\geqslant0}$.

\begin{notation}\label{not:dmax}
For an integer $m\geqslant1$, put
\[
 d_{\max}:=m-1-(q-1)\left\lceil\frac{m}{q}\right\rceil .
\]
\end{notation}

\begin{lemma}\label{lem:semigroup}
\textup{(i)} The set $\Nq$ has the explicit interval description
\[
  \Nq=\bigcup_{j\geqslant0}[\,2q+j(q-1),2q+jq\,]\cap\mathbb Z.
\]
Moreover, for every $m\geqslant1$, $d_{\max}>0$ if and only if
$m\in\Nq$.
 For an integer $m\geqslant2q$, write $m-2q=aq+r$ with $a\geqslant0$
and $0\leqslant r<q$. Then $m\notin\Nq$ exactly when $1\leqslant r$ and
$a+r\leqslant q-2$. 
\par\medskip\noindent\textup{(ii)}  Every $m\geqslant q^2-q+2$ belongs to $\Nq$. The complement
among integers $m\geqslant2q$ is the explicit finite list  
 \[ \{m\in\mathbb Z:m\geqslant2q,\ m\notin\Nq\}
 =\{2q+aq+r:0\leqslant a\leqslant q-3,\ 1\leqslant r\leqslant q-2-a\},
\]
 of cardinality $(q-1)(q-2)/2$; it is empty for $q=2$. 

\end{lemma}

\begin{proof}

\emph{(i).} If $m=2q+aq+b(q-1)$ with $a,b\geqslant0$, put $j:=a+b$. Then
$m=2q+jq-b$ and $0\leqslant b\leqslant j$, which gives the displayed interval.
Conversely, if $2q+j(q-1)\leqslant m\leqslant2q+jq$, write $m=2q+jq-b$ with
$0\leqslant b\leqslant j$; then $m=2q+(j-b)q+b(q-1)$.

We now prove the sign criterion. Put $u:=\lceil m/q\rceil$.
Since $m\leqslant uq$, the condition $d_{\max}>0$ is equivalent to
$(q-1)u+2\leqslant m\leqslant uq$. These inequalities are impossible
for $u=1$. For $u\geqslant2$, setting $j:=u-2$ identifies them with
the displayed interval, so $d_{\max}>0$ implies $m\in\Nq$.
Conversely, if $m$ lies in the interval indexed by some $j\geqslant0$,
then $m\leqslant(j+2)q$, hence $u\leqslant j+2$, and
\[
 m\geqslant2q+j(q-1)=(q-1)(j+2)+2
 \geqslant(q-1)u+2,
\]
so $d_{\max}>0$. This proves the asserted equivalence.

For the remainder criterion, put $n:=m-2q=aq+r$.
If $r=0$, then $n=aq$ is a sum of the two generators.
If $1\leqslant r<q$, a representation $n=uq+v(q-1)$ with
$u,v\in\mathbb Z_{\geqslant0}$ requires $v\equiv-r\pmod q$, hence $v\geqslant q-r$.
Such a representation exists if and only if
$n\geqslant(q-r)(q-1)$: the difference $n-(q-r)(q-1)$ is divisible by $q$.
Substituting $n=aq+r$ rewrites this inequality as $a+r\geqslant q-1$.

\emph{(ii).} For $j\geqslant q-2$, consecutive intervals satisfy
$(j+1)(q-1)\leqslant jq+1$, so their union covers every
$m\geqslant2q+(q-2)(q-1)=q^2-q+2$. By the remainder criterion, $m\notin\Nq$
requires $a\leqslant q-3$ and $r\leqslant q-2-a$, and each such pair gives one
exception. The count is
$\sum_{a=0}^{q-3}(q-2-a)=(q-1)(q-2)/2$. For $q=2$, the list is empty.
\end{proof}

\begin{samepage}
\begin{proposition}\label{prop:trichotomy}
Let $S$ be nonempty, put $m:=|S|$, and let $d_{\max}$ be as in
Notation~\ref{not:dmax}. Then:
\begin{enumerate}[label=(\roman*)]
\item The condition $d_{\max}<0$ holds if and only if $G_S=0$,
equivalently, $\mathbb{E}_{S,d}=\varnothing$ for every
$d\geqslant0$.
\item The condition $d_{\max}>0$ holds if and only if $G_S$ is
nonconstant, equivalently, $m\in\Nq$. If $G_S\ne0$, then
$\deg_X G_S=d_{\max}$.
\item The condition $d_{\max}=0$ holds if and only if $G_S$ is a
nonzero constant. This occurs exactly when
$m=q+(q-1)k$ with $0\leqslant k<q$, equivalently,
$\mathbb{E}_{S,0}\ne\varnothing$ and $m\notin\Nq$.
\end{enumerate}
\end{proposition}
\end{samepage}

\begin{proof}

By Notation~\ref{not:scalar-to-differential},
$\mathfrak g_{S,d}=g_{S,d}\nu_S$ with $\nu_S\ne0$. Thus
Lemma~\ref{lem:support} also characterizes the supported degrees of
$G_S$.

\emph{(i).} For fixed $m$, let $r\in\{0,\ldots,q-1\}$ satisfy
$r\equiv-m\pmod q$.
The identity $r=q\lceil m/q\rceil-m$ shows that the
definition in Notation~\ref{not:dmax} is equivalent to
$d_{\max}=(m-(q-1)r)/q-1$.
The support condition in Lemma~\ref{lem:support} requires
$k=r+qj$ with $j\geqslant0$, and the corresponding degree is
$d=d_{\max}-(q-1)j$. Thus no nonnegative degree exists exactly when
$d_{\max}<0$, and then $G_S=0$. If $d_{\max}\geqslant0$, a
nonzero coefficient exists, so $G_S\ne0$. The equivalence with
$\mathbb{E}_{S,d}=\varnothing$ for every $d\geqslant0$ follows
from Theorem~\ref{thm:formula} and Lemma~\ref{lem:depth-existence}.

\emph{(ii).} For $G_S\ne0$, the computation in part~(i)
shows that its largest supported degree is $d_{\max}$. Thus
$\deg_XG_S=d_{\max}$, and $G_S$ is nonconstant if and only if
$d_{\max}>0$. Lemma~\ref{lem:semigroup}(i) identifies this last
condition with $m\in\Nq$.

\emph{(iii).} By (i) and (ii), $G_S$ is a nonzero constant
exactly in the remaining case $d_{\max}=0$.
Degree zero is supported exactly when $m=q+(q-1)k$ for some
$k\geqslant0$. Fix such a $k$. Any supported degree $d$, with associated
integer $k'\geqslant0$, satisfies $qd=(q-1)(k-k')$. Since $q$ and $q-1$
are coprime, this is equivalent to $d=(q-1)j$ and $k'=k-qj$ for an
integer $j\geqslant0$. A positive supported degree occurs exactly when
$k\geqslant q$: necessity follows from $k-qj=k'\geqslant0$, and
sufficiency follows by taking $j=1$. Thus $G_S$ is a nonzero constant
exactly when $0\leqslant k<q$. These conditions say precisely that
$\mathbb{E}_{S,0}\ne\varnothing$ and $m\notin\Nq$.
\end{proof}

\begin{remark}
Every $m\geqslant q^2-q+2$ gives a nonconstant correction,
while the constant cases with $m\geqslant2q$ are precisely the finite
list in Lemma~\ref{lem:semigroup}(ii).
\end{remark}

\begin{corollary}\label{cor:simultaneous-constancy}
All canonical corrections $\mathfrak{G}_S$ for nonempty
$S\subseteq\Sigma_s$ are constant if and only if $s<2q$.
\end{corollary}

\begin{proof}
 If $s<2q$, every nonempty subset has $m<2q$, so $m\notin\Nq$,
and Lemma~\ref{lem:support} makes its correction constant. If
$s\geqslant2q$, a subset of size $2q$ has $2q\in\Nq$, and
Proposition~\ref{prop:trichotomy}(ii) gives a nonzero degree-one
coefficient. 
\end{proof}

\begin{corollary}\label{cor:degrees}
Let $d_{\max}$ be as in Notation~\ref{not:dmax},
with $d_{\max}\geqslant0$. The leading coefficient of $\mathfrak{G}_{S}$
has the following explicit forms:
\begin{align}
 \mathfrak{g}_{S,d_{\max}}
   &=\prodext_{i\in S}\Delta_{\varphi,i}
       &&(m\equiv0\pmod q),\label{eq:leading}\\
 \mathfrak{g}_{S,d_{\max}}&=
 \sum_{\substack{T\subseteq S\\|T|=q}}
    \left(\prodext_{i\in T}\mathfrak{B}_{i}^{[1]}\right)
    \left(\prodext_{i\in S\setminus T}\Delta_{\varphi,i}\right)
       &&(m\equiv-1\pmod q).\label{eq:second-leading}
\end{align}
\end{corollary}

\begin{proof}

In the first case, $m=qr$ for some $r\geqslant1$, so
$\lceil m/q\rceil=r$ and $d_{\max}=r-1$; hence
$m=q(d_{\max}+1)$. The weight equation is
$\sum_{i\in S}q^{-n_i}=m$. Since $q^{-n_i}\leqslant1$ for every $i$, all
terms must equal $1$, so every $n_i=0$. This gives \eqref{eq:leading}.

In the second case, $m=qr+q-1$ for some $r\geqslant0$, so
$\lceil m/q\rceil=r+1$ and $d_{\max}=r-1$; hence
$m=q(d_{\max}+1)+q-1$. The depth bound in
Lemma~\ref{lem:chain-contribution}(ii) gives $n_i\in\{0,1\}$.
Put $T:=\{i:n_i=1\}$. The weight equation becomes
\[
 m-(1-q^{-1})|T|=q(d_{\max}+1)=m-(q-1),
\]
so $|T|=q$. Summing the corresponding products in
\eqref{eq:coefficient-formula} gives \eqref{eq:second-leading}.
\end{proof}

\subsection{Definition field of the correction polynomials}

Corollary~\ref{cor:rational-coefficients} gives rational correction
coefficients over $\Cinf$. We descend them to the module's definition
field by computing their residue-dual coordinates as finite torsion traces.

The continuous
residue pairing of \cite[Section 3.1]{HYY} satisfies
\[
 \langle u,\Hcal(\lambda)\rangle_\infty
       =e(\iota(u)\lambda)\qquad(u\in K_\infty,\ \lambda\in\Lambda).
\]
For $\eta\in\Omega_A$, the pairing is
$\Tr_{\mathbb F_\infty/\Fq}\Res_\infty(u\eta)$, extended
$\Cinf$-linearly and continuously to $\mathbb M_A$.
 
Taking its product over the factors in a nonempty factor set $S$ gives a continuous pairing
on $\mathbb M_S$. On pure tensors,
\[
 \left\langle(u_i)_{i\in S},\prodext_{i\in S}\eta_i\right\rangle
 =\prod_{i\in S}\langle u_i,\eta_i\rangle_\infty
       \qquad(u_i\in K_\infty,\ \eta_i\in\mathbb M_A).
\]
Continuity permits termwise pairing in the normally convergent
sums below.
 
\begin{lemma}\label{lem:finite-torsion-pairing}
Suppose that $\varphi:A\to H\{\tau\}$ is defined over $H$, where
$\iota(K)\subset H\subset\Cinf$. Fix $b\in A\setminus\Fq$ and put
$\beta:=\iota(b)$. For $c\in A$, let $J_{b,c}(Z)\in H[Z]$ be the
remainder of $\varphi_c(Z)$ modulo $\varphi_b(Z)$. For every nonempty
$S$ and all $r_i\in A$ ($i\in S$), coefficientwise pairing gives
\begin{equation}\label{eq:finite-descent-pairing}
 \left\langle(r_i/b)_{i\in S},\mathfrak G_S(X)\right\rangle
 =-\frac1\beta\Tr_{\mathcal A_{b,X}/H[X]}
       \left(\frac{\prod_{i\in S}J_{b,r_i}(Z)}Z\right),
\end{equation}
where
\begin{equation}
 \mathcal A_{b,X}:=H[X,Z]/(\varphi_b(Z)+X).
\end{equation}
The trace argument is a polynomial, and the right-hand side belongs to
$H[X]$.
\end{lemma}

\begin{proof}
The polynomial $J_{b,c}$ is a $q$-polynomial with zero constant term:
division proceeds by subtracting scalar multiples of
$\varphi_b(Z)^{q^j}$. Fix $w\notin\Lambda$ and put
$t:=e(w/\beta)$, so that $\varphi_b(t)=e(w)$. The map $\lambda\mapsto e(\lambda/\beta)$
identifies $\Lambda/\beta\Lambda$ with the torsion space
$\varphi[b]:=\{x:\varphi_b(x)=0\}$. The map from $\Lambda$ has kernel
$\beta\Lambda$. For $x\in\varphi[b]$, surjectivity of $e$ supplies
$u\in\Cinf$ with $e(u)=x$;
then $e(\beta u)=\varphi_b(x)=0$, so $\beta u\in\Lambda$.
For a representative $\lambda_0$ and $x=e(\lambda_0/\beta)$,
the scalar partial fraction identity gives
\[
 \sum_{\mu\in\beta\Lambda}\frac1{w-\lambda_0-\mu}
 =\frac1{\beta e((w-\lambda_0)/\beta)}
 =\frac1{\beta(t-x)}.
\]
Pairing \eqref{eq:inputs} with $c/b$ and grouping by these cosets gives
\[
 \langle c/b,\Hcal(w)\rangle_\infty
 =\frac{\varphi_b(t)}\beta
       \sum_{x\in\varphi[b]}\frac{\varphi_c(x)}{t-x}
 =J_{b,c}(t).
\]
For $\lambda$ in the coset corresponding to $x$, residue
duality and the exponential equation give
$\langle c/b,\Hcal(\lambda)\rangle_\infty
=e(\iota(c)\lambda/\beta)=\varphi_c(x)$.
The last equality in the display follows from
Lemma~\ref{lem:interpolation-values}: $\varphi_b/\beta$ is the
normalized subspace polynomial of $\varphi[b]$, and $J_{b,c}$ is the
unique polynomial of degree less than $\deg\varphi_b$ that agrees
with $\varphi_c$ on this space.

For $r_i\in A$, pairing \eqref{eq:canonical} with $(r_i/b)_{i\in S}$
now gives numerator $\prod_{i\in S}J_{b,r_i}(x-t)$ on the coset
corresponding to $x$. Its denominator sum is $1/(\beta(t-x))$.
The elements $x-t$ are the roots of $\varphi_b(Z)+e(w)$. Consequently,
evaluating at $X=e(w)$ and summing over $x$ gives
\eqref{eq:finite-descent-pairing}.

Each $J_{b,r_i}$ is divisible by $Z$, and $S$ is nonempty, so the trace
argument is a polynomial. The algebra is finite free because the
leading coefficient of $\varphi_b$ is a unit of $H$, and it is \'etale
because $\varphi_b'=\beta\ne0$. Thus the right side belongs to $H[X]$.
The calculation proves equality at $X=e(w)$ for $w\notin\Lambda$;
these values form an infinite set, so it proves the polynomial identity.
\end{proof}

For $0\ne b\in A$, put $\deg(bA):=\dim_{\Fq}(A/bA)$.
The following residue-duality result supplies a common denominator for
the coordinate functionals. Duals are taken over $\Fq$.
\begin{lemma}[{\cite[Lemma 3.3]{HYY}}]\label{lem:hyy-residue-surjectivity}
Let $D=M\infty$ with an integer $M\geqslant1$. If $b\in A\setminus\{0\}$ satisfies
\[
 \deg(bA)>2g-2+MN,
\]
then the residue map
\[
 A/bA\longrightarrow H^0({\mathcal X},\Omega_{\mathcal X}^1(D))^*,\qquad
 c+bA\longmapsto\bigl(\nu\longmapsto
 \langle c/b,\nu\rangle_\infty\bigr)
\]
is surjective.
\end{lemma}

\begin{proposition}\label{prop:correction-field-rationality}
Suppose that $\varphi:A\to H\{\tau\}$ is defined over $H$, where
$\iota(K)\subset H\subset\Cinf$. For every nonempty $S$, every
coefficient of $\mathfrak{G}_S$ is a rational external differential over $H$.
\end{proposition}

\begin{proof}
Fix a nonzero coefficient $\mathfrak g_{S,d}$.
Lemma~\ref{lem:support} gives a nonnegative integer
$k=(|S|-q(d+1))/(q-1)$. Each tuple $(n_i)\in\mathbb E_{S,d}$
in its expansion satisfies $n_i\leqslant k$ by
Lemma~\ref{lem:chain-contribution}(ii).
By \eqref{eq:W} and \eqref{eq:Gamma-divisor},
\[
 (\mathfrak{B}^{[n]})=V_\varphi+W^{(-n)}+\sum_{j=1}^n\Xi^{(-j)}
                         -\sum_{j=1}^{n+2}P^{(-j)}.
\]
Each point over $\infty$ occurs at most
$\lceil(n+2)/N\rceil$ times in the negative sum. Thus, with
\[
 M:=\lceil(k+2)/N\rceil,\qquad
 E_M:=H^0({\mathcal X},\Omega_{\mathcal X}^1(M\infty)),
\]
every basis differential occurring in $\mathfrak g_{S,d}$ belongs
to $H^0({\mathcal X}_{\Cinf},\Omega^1(M\infty))$.
Base change identifies this finite-dimensional space with
$\Cinf\otimes_{\Fq}E_M$. Thus $\mathfrak g_{S,d}$ lies in the
external tensor product of these spaces over the factors in $S$,
with tensors taken over $\Cinf$.

Choose an $\Fq$-basis $\nu_1,\ldots,\nu_\ell$ of $E_M$.
Riemann--Roch supplies $b\in A\setminus\Fq$ with
$\deg(bA)>2g-2+MN$. Set $\beta:=\iota(b)$, and use
$J_{b,c}$ and $\mathcal A_{b,X}$ from
Lemma~\ref{lem:finite-torsion-pairing}.
Lemma~\ref{lem:hyy-residue-surjectivity} supplies $r_1,\ldots,r_\ell\in A$
representing the dual basis:
$\langle r_j/b,\nu_l\rangle_\infty=\delta_{jl}$.
Write
$\mathfrak g_{S,d}=\sum_{1\leqslant j_i\leqslant\ell}
\alpha_{(j_i)}\prodext_{i\in S}\nu_{j_i}$,
where $\alpha_{(j_i)}\in\Cinf$.
Pairing with $(r_{j_i}/b)_{i\in S}$ extracts $\alpha_{(j_i)}$.
By \eqref{eq:finite-descent-pairing}, $\alpha_{(j_i)}$ is the
coefficient of $X^d$ in
\begin{equation}\label{eq:paper1-descent-coordinates}
 -\frac1\beta\Tr_{\mathcal A_{b,X}/H[X]}
       \left(\frac{\prod_{i\in S}J_{b,r_{j_i}}(Z)}Z\right)
       \in H[X].
\end{equation}
Thus every $\alpha_{(j_i)}$ belongs to $H$, proving the descent.
\end{proof}

\section{Multicharacter Pellarin series and class-field trace identities}\label{sec:mc-traces}

To obtain the trace over the Hilbert class field, we descend the
multicharacter Pellarin values to the module's definition field and
compare the lattice sums with the ideal series.

Let $\widehat\Dcal_S$ be the completion of $\Dcal_S$ for the
Gauss norm transported by $\eta\mapsto\eta/\nu_S$.
The comparison following Notation~\ref{not:scalar-trivialization}
identifies $\mathbb M_S$ with a closed subspace of $\widehat\Dcal_S$.

\subsection{The \texorpdfstring{$L$}{L}-value expansion}\label{sec:L-expansion}

\begin{definition}\label{def:mc-pellarin}
For every $n\geqslant1$, the multicharacter Ferraro series of
weight $n$ in differential normalization is
\[
 C_{\varphi,n}^{s}:=\sum_{0\ne\lambda\in\Lambda}
       \frac{\Hcal_1(\lambda)\boxtimes\cdots\boxtimes\Hcal_s(\lambda)}
            {\lambda^n}.
\]
\end{definition}

\Needspace{9\baselineskip}
\begin{proposition}\label{prop:L-expansion}
\textup{(i)} The series $C_{\varphi,n}^{s}$ converges in
$\widehat\Dcal_{\Sigma_s}$ for every $n\geqslant1$, and $C_{\varphi,n}^{s}=0$ unless
$n\equiv s\pmod{q-1}$.
\par\medskip\noindent\textup{(ii)} For $|w|<\min_{0\ne\lambda\in\Lambda}|\lambda|$,
\[
 \mathscr C_{\varphi}^s(w)
 =-\sum_{\substack{n\geqslant1\\n\equiv s\ (q-1)}}C_{\varphi,n}^{s}w^{n-1}.
\]

\end{proposition}
\begin{proof}

\emph{(i).} By Lemma~\ref{lem:period-bound}(ii), the scalar products
$\prod_i h_i(\lambda)$ have a common Gauss-norm bound $M$.
After division by $\nu_{\Sigma_s}$, the summand of $C_{\varphi,n}^{s}$
therefore has norm at most $M/|\lambda|^n$. Since $n\geqslant1$,
these norms tend to zero along the discrete lattice, proving convergence.
For $\zeta\in\Fq^\times$, Proposition~\ref{prop:fractional-linearity}(i) gives
$\Hcal_i(\zeta\lambda)=\zeta\Hcal_i(\lambda)$, and therefore
\[
 \frac{\prodext_{i=1}^s\Hcal_i(\zeta\lambda)}{(\zeta\lambda)^n}
 =\zeta^{s-n}\frac{\prodext_{i=1}^s\Hcal_i(\lambda)}{\lambda^n}.
\]
Reindexing the convergent sum by $\lambda\mapsto\zeta\lambda$ gives
$C_{\varphi,n}^{s}=\zeta^{s-n}C_{\varphi,n}^{s}$. If
$n\not\equiv s\pmod{q-1}$, a generator of $\Fq^\times$ has
$\zeta^{s-n}\ne1$, which forces $C_{\varphi,n}^{s}=0$.

\emph{(ii).} Fix $0<\rho<r_\varphi$, where
$r_\varphi=\min_{0\ne\lambda\in\Lambda}|\lambda|$.
Retain the bound $M$ from part~(i) and use the disk Gauss norm
$\norm{\cdot}_{F,\rho}$ from Section~\ref{sec:limits}.
The terms of the geometric expansion satisfy
\[
 \left\lVert
  \frac{\prod_{i=1}^s h_i(\lambda)}{\lambda}
             \left(\frac w\lambda\right)^j
 \right\rVert_{F,\rho}
 \leqslant\frac{M}{|\lambda|}
                \left(\frac\rho{|\lambda|}\right)^j.
\]
For each $\varepsilon>0$, the condition that the bound exceeds
$\varepsilon$ bounds $|\lambda|$. Since $\rho/r_\varphi<1$, it also
bounds $j$. Thus only finitely many pairs $(\lambda,j)$ contribute
terms of norm greater than $\varepsilon$, and the double family is
summable in the complete disk algebra.
We may therefore exchange the sums and apply part~(i):
\[
\begin{aligned}
 \mathscr C_\varphi^s(w)
 &=\nu_{\Sigma_s}\sum_{0\ne\lambda\in\Lambda}
   \frac{\prod_{i=1}^s h_i(\lambda)}{w-\lambda}\\
 &=-\nu_{\Sigma_s}\sum_{0\ne\lambda\in\Lambda}\sum_{j\geqslant0}
   \frac{\prod_{i=1}^s h_i(\lambda)}{\lambda}
   \left(\frac w\lambda\right)^j\\
 &=-\sum_{n\geqslant1}C_{\varphi,n}^{s}w^{n-1}\\
 &=-\sum_{\substack{n\geqslant1\\n\equiv s\ (q-1)}}
   C_{\varphi,n}^{s}w^{n-1}.
\end{aligned}
\]
\end{proof}

The following projection extracts the Taylor terms with
$n\equiv s\pmod{q^N-1}$ and will be used in the Hilbert-class trace
formula of Theorem~\ref{thm:mc-traceE}.
\begin{notation}\label{not:mc-projection}
Define the coefficient projector by
\[
 (\mathcal P_{q^N}^{s}B)(w)
 =\frac1{q^N-1}\sum_{\varepsilon\in\iota(\mathbb F_\infty^\times)}
            \varepsilon^{1-s}B(\varepsilon w).
\]
\end{notation}
\begin{proposition}\label{thm:mc-trace}
Near $w=0$, we have the expansion
\begin{equation}\label{eq:mc-projection}
 (\mathcal P_{q^N}^{s}\mathscr C_{\varphi}^s)(w)
 =-\sum_{\substack{n\geqslant1\\n\equiv s\ (q^N-1)}}
          C_{\varphi,n}^{s}w^{n-1}.
\end{equation}
\end{proposition}
\begin{proof}

Proposition~\ref{prop:L-expansion}(ii) and the
definition of $\mathcal P_{q^N}^{s}$ give
\[
 (\mathcal P_{q^N}^{s}\mathscr C_{\varphi}^s)(w)
 =-\sum_{\substack{n\geqslant1\\n\equiv s\ (q-1)}}
 C_{\varphi,n}^{s}w^{n-1}
 \left(\frac1{q^N-1}\sum_{\varepsilon\in\iota(\mathbb F_\infty^\times)}
 \varepsilon^{n-s}\right).
\]
The inner character average is
\[
 \frac1{q^N-1}\sum_{\varepsilon\in\iota(\mathbb F_\infty^\times)}
 \varepsilon^{n-s}
 =
 \begin{cases}
  1,&(q^N-1)\mid(n-s),\\
  0,&(q^N-1)\nmid(n-s).
 \end{cases}
\]
When $(q^N-1)\nmid(n-s)$, multiplication by an
$\varepsilon_0$ with $\varepsilon_0^{n-s}\ne1$ forces the average to
vanish; otherwise all its summands are one.
Substituting this average gives \eqref{eq:mc-projection}.
\end{proof}

\subsection{Comparison with element series}\label{sec:mc-element-comparison}
Fix a sign function
$\operatorname{sgn}:K_\infty^\times\to\mathbb F_\infty^\times$.

\begin{definition}\label{def:ele-L-series}
For a fractional ideal $\mathfrak j$, write
$\mathfrak j_+:=\{c\in\mathfrak j:\operatorname{sgn}(c)=1\}$ and define
\[
 L^s(\mathfrak j;n)
 :=\sum_{0\ne c\in\mathfrak j}
       \frac{\prod_i\chi_i(c)}{\iota(c)^n},\qquad
 L^s(\mathfrak j_+;n)
 :=\sum_{c\in\mathfrak j_+}
       \frac{\prod_i\chi_i(c)}{\iota(c)^n}.
\]
Both sums converge in the scalar Gauss completion: clearing a fixed
denominator of $\mathfrak j$ bounds the geometric numerators uniformly.
\end{definition}
\begin{proposition}\label{prop:mc-element-normalization}
Suppose that $\Lambda_\psi=\pC\iota(A)$, and put
$\mathscr W_{\psi,s}:=\prodext_i\Hcal_{\psi,i}(\pC)$.
\begin{enumerate}[label=(\roman*)]
\item For every $n\geqslant1$,
\begin{equation}\label{eq:mc-element-normalization-all}
 C_{\psi,n}^{s}=\pC^{-n}\mathscr W_{\psi,s}\,
 {L^s}(A;n).
\end{equation}
\item If $N=1$ and $n\equiv s\pmod{q-1}$, then
\begin{equation}\label{eq:mc-element-normalization}
 C_{\psi,n}^{s}=-\pC^{-n}\mathscr W_{\psi,s}\,{L^s}(A_+;n).
\end{equation}
\end{enumerate}
\end{proposition}
\begin{proof}

\emph{(i).} Every nonzero period has the form
$\lambda=\pC\iota(a)$ with $a\in A\setminus\{0\}$, and
\eqref{eq:fractional-period-linearity} gives
$\Hcal_{\psi,i}(\pC\iota(a))=\chi_i(a)\Hcal_{\psi,i}(\pC)$. Consequently,
\[
\begin{aligned}
 C_{\psi,n}^{s}
 &=\sum_{0\ne a\in A}
 \frac{\prodext_{i=1}^s\Hcal_{\psi,i}(\pC\iota(a))}
      {\pC^n\iota(a)^n}\\
 &=\pC^{-n}\mathscr W_{\psi,s}
 \sum_{0\ne a\in A}\frac{\prod_i\chi_i(a)}{\iota(a)^n}
 =\pC^{-n}\mathscr W_{\psi,s}\,L^s(A;n),
\end{aligned}
\]
which is \eqref{eq:mc-element-normalization-all}.

\emph{(ii).} Because $N=1$, each $\Fq^\times$-orbit in
$A\setminus\{0\}$ contains exactly one sign-one representative $a_0$.
Grouping the full element sum into these orbits gives
\[
\begin{aligned}
 L^s(A;n)
 &=\sum_{a_0\in A_+}\sum_{\varepsilon\in\Fq^\times}
 \frac{\prod_i\chi_i(\varepsilon a_0)}
      {\iota(\varepsilon a_0)^n}\\
 &=\sum_{a_0\in A_+}
 \left(\sum_{\varepsilon\in\Fq^\times}\varepsilon^{s-n}\right)
 \frac{\prod_i\chi_i(a_0)}{\iota(a_0)^n}\\
 &=(q-1)L^s(A_+;n)=-L^s(A_+;n),
\end{aligned}
\]
where the last two equalities use $n\equiv s\pmod{q-1}$ and $q-1=-1$
in characteristic $p$. Part~(i) now gives
\eqref{eq:mc-element-normalization}.
\end{proof}

The following Carlitz normalization identifies these lattice
values with the multicharacter values of Pellarin--Perkins.
\begin{notation}\label{not:carlitz-normalization}
For the Carlitz module over $A:=\Fq[\theta]$, let $\exp_C$ be
the Carlitz exponential and choose $\pC\in\Cinf^\times$ generating its
period lattice. For $k\geqslant0$, put
\[
 e_k(z):=\exp_C(\pC z/\theta^k),
\]
so $e_C:=e_0$ satisfies $\ker e_C=A$. Anderson's generating function is
\[
 f_t(z):=\sum_{k\geqslant1}t^{k-1}e_k(z).
\]
This series converges normally on closed $z$-disks whenever
$|t|\leqslant1$; see \cite[Section 1.1]{PP16}. We use the normalizations
\[
 \Hcal(w)=-f_t(w/\pC)\,dt,\qquad
 C_\varphi=\frac{dt}{t-\theta},\qquad \Delta_\varphi=dt.
\]
Choose $\nu:=dt$, so $h_i(w)=-f_{t_i}(w/\pC)$.
\end{notation}

\begin{remark}\label{rem:L-identity}
In the Carlitz normalization of
Notation~\ref{not:carlitz-normalization}, the expansion in
Proposition~\ref{prop:L-expansion}(ii) recovers
\cite[Lemma 20]{PP16}.

The multicharacter $L$-function of \cite{PP16} is
\[
 L(\chi_{t_1},\ldots,\chi_{t_s};n)
 :=\sum_{a\in A_+}\frac{\chi_{t_1}(a)\cdots\chi_{t_s}(a)}{a^n},
\]
where $
 A_+=\{a\in\Fq[\theta]\ \mbox{\normalfont monic}\}$. For Anderson's special function $\omega$, one has
$\chi_{t_i}(a)=\omega(t_i)^{-1}f_{t_i}(a)$ \cite{PP16}. At
$n\equiv s\pmod{q-1}$, the substitutions $\lambda=\pC a$ and
$\Hcal_i(\pC a)=-f_{t_i}(a)\,dt_i=-\omega(t_i)\chi_{t_i}(a)\,dt_i$
give
\[
 C_{\varphi,n}^{s}=\frac{(-1)^{s+1}\omega(t_1)\cdots\omega(t_s)}{\pC^{\,n}}\,
       L(\chi_{t_1},\ldots,\chi_{t_s};n)\prodext_{i=1}^s dt_i.
\]

For an arbitrary function ring $A$ and $s=1$, the correction
$\mathfrak{G}_{\{1\}}$ vanishes by the
degree bound, and the identity reduces to the partial fraction
expansion \eqref{eq:inputs}:
\[
 \frac{\Hcal(w)}{e(w)}
 =-\sum_{\substack{n\geqslant1\\n\equiv1\ (q-1)}}C_{\varphi,n}^{1}w^{n-1}.
\]
This recovers \cite[Theorem D]{HYY}. Taking $w\to0$ and using
$e'(0)=1$ and $\Hcal'(0)=C_\varphi$ gives
$C_\varphi=-C_{\varphi,1}^{1}$.
\end{remark}

\subsection{Field covariance of the multicharacter Pellarin series} \label{sec:field-covariance}
\leavevmode
Let $\varphi$ be defined over a field $H$ with
$\iota(K)\subset H\subset\Cinf$. We use the generating identity to
descend the lattice values to $H$ and to compare them under embeddings
of $H$.
\begin{lemma}\label{lem:exponential-field-rationality}
Write $e(w)=\sum_{r\geqslant0}e_rw^{q^r}$, with $e_0=1$. Then $e_r\in H$
for all $r\geqslant0$.
\end{lemma}
\begin{proof}
For any
$a_0\in A\setminus\Fq$, write
$\varphi_{a_0}(Z)=\sum_{\ell=0}^{d_0}\kappa_\ell Z^{q^\ell}$,
where $d_0:=\deg_\tau\varphi_{a_0}$, $\kappa_0=\iota(a_0)$, and
$\kappa_\ell\in H$. Comparing the coefficient of $w^{q^r}$ in
$e(\iota(a_0)w)=\varphi_{a_0}(e(w))$ gives
\[
 \bigl(\iota(a_0)^{q^r}-\iota(a_0)\bigr)e_r
 =\sum_{\ell=1}^{\min\{d_0,r\}}
 \kappa_\ell e_{r-\ell}^{q^\ell}
 \qquad(r\geqslant1).
\]
Since $a_0$ is transcendental over $\Fq$, the factor
$\iota(a_0)^{q^r}-\iota(a_0)$ is nonzero. Induction from $e_0=1$
therefore gives $e_r\in H$ for every $r\geqslant0$.
\end{proof}

\begin{samepage}
\begin{proposition}\label{prop:field-covariance}
Suppose that $\varphi:A\to H\{\tau\}$ is defined over $H$, where
$\iota(K)\subset H\subset\Cinf$.
\begin{enumerate}[label=(\roman*)]
\item Every $C_{\varphi,n}^{s}$ is a rational external differential over $H$.
\item For an $\iota(K)$-embedding $\sigma:H\hookrightarrow\Cinf$,
let $\varphi^\sigma$ be obtained by applying $\sigma$ to the
coefficients of $\varphi$. Then
\[
 \sigma(C_{\varphi,n}^{s})=C_{\varphi^\sigma,n}^{s}.
\]
\item For $a\in\Cinf^\times$ and $\varphi^a=a\varphi a^{-1}$,
\begin{equation}\label{eq:mc-scaling}
 C_{\varphi^a,n}^{s}=a^{s-n}C_{\varphi,n}^{s}.
\end{equation}
\end{enumerate}
\end{proposition}
\end{samepage}

\begin{proof}
\emph{(i).}
Write $\Hcal(w)=\sum_{r\geqslant0}T_rw^{q^r}$, and retain the
exponential coefficients $e_r\in H$ from
Lemma~\ref{lem:exponential-field-rationality}.
By \cite[Corollary 3.8 and Lemma 3.9]{HYY}, $f$ and $C_\varphi$ are
rational over $H$. The one-factor recursion
\[
 T_0=C_\varphi,\qquad
 T_r=e_rC_\varphi+f^{-1}T_{r-1}^{(1)}\quad(r\geqslant1)
\]
from \cite[Proposition 4.8]{HYY} gives
$T_r\in\Omega^1_{H({\mathcal X})/H}$ by induction, since coefficient
Frobenius preserves $H$.

To extract Taylor coefficients from \eqref{eq:generating}, write
$e(w)=wE(w)$ with $E(w)\in1+w^{q-1}H[[w]]$.
Recursive inversion gives $E(w)^{-1}\in H[[w]]$. The term indexed by
$S=\varnothing$ becomes
\[
 w^{s-1}E(w)^{-1}\prodext_{i=1}^s\frac{\Hcal_i(w)}w,
\]
which is regular at zero. Every other term is a polynomial in $e(w)$
times generating differentials. Each Taylor coefficient is therefore a finite expression in
the $e_r$, the rational differentials $T_r$, and the correction
coefficients. The latter are rational over $H$ by
Proposition~\ref{prop:correction-field-rationality}. Thus all Taylor
coefficients are rational over $H$, and
Proposition~\ref{prop:L-expansion}(ii) identifies them with
$-C_{\varphi,n}^{s}$.

\emph{(ii).} Extend $\sigma$ by fixing geometric functions
and differentials. For each nonzero correction coefficient indexed by
$(S,d)$, Proposition~\ref{prop:correction-field-rationality} uses the
pole bound $M=\lceil(k+2)/N\rceil$, where
$k=(|S|-q(d+1))/(q-1)$. This bound depends only on $|S|$, $d$, and $N$.
We may therefore use the same space $E_M$, $\Fq$-basis, and residue-dual
classes in \eqref{eq:paper1-descent-coordinates} for $\varphi$ and
$\varphi^\sigma$. The scalar $\beta=\iota(b)$ is also fixed by $\sigma$.
Applying $\sigma$ to the polynomial remainders and finite free traces in
\eqref{eq:paper1-descent-coordinates} therefore shows that
$(\mathfrak G_S)^\sigma$, obtained by applying $\sigma$ coefficientwise,
is the correction polynomial of $\varphi^\sigma$.
Lemma 3.9 of \cite{HYY} also gives
$f_{\varphi^\sigma}=\sigma(f)$ and
$C_{\varphi^\sigma}=\sigma(C_\varphi)$.
The exponential recursion in
Lemma~\ref{lem:exponential-field-rationality} and the differential
recursion in part~(i) then identify the Taylor coefficients of
$e_{\varphi^\sigma}$ and $\Hcal_{\varphi^\sigma}$ with the
$\sigma$-transforms of those of $e$ and $\Hcal$.
Recursive inversion of $e(w)/w$ does the same for the empty-subset term
in \eqref{eq:generating}.

Applying $\sigma$ coefficientwise to the Taylor expansion of
\eqref{eq:generating} now gives
\[
\begin{aligned}
-\sum_{n\geqslant1}\sigma(C_{\varphi,n}^{s})w^{n-1}
&=\sum_{S\subseteq\Sigma_s}
(\mathfrak G_S)^\sigma(e_{\varphi^\sigma}(w))
\boxtimes\prodext_{i\notin S}\mathcal H_{\varphi^\sigma,i}(w)\\
&=\mathscr C_{\varphi^\sigma}^s(w)
=-\sum_{n\geqslant1}C_{\varphi^\sigma,n}^{s}w^{n-1}.
\end{aligned}
\]
Comparing the coefficients of $w^{n-1}$ proves the asserted covariance.

\emph{(iii).} The lattice of $\varphi^a$ is $a\Lambda$, and
$\Hcal_{\varphi^a}(a\lambda)=a\Hcal_\varphi(\lambda)$ by
\cite[equation (3.4)]{HYY}. Substitution in the defining period sum gives
\eqref{eq:mc-scaling}.
\end{proof}

\subsection{Ideal series and class decomposition}\label{sec:mc-trace}

Fix compatible embeddings of the ordinary and sign Hilbert class
fields $H_A\subset H_A^+$ into $\Cinf$ extending $\iota|_K$,
and identify them with their images. The notation $\Tr_{H_A/K}$ denotes the
trace over $\iota(K)$, with the geometric variables fixed.

Fix a Hayes module $\psi$ over the sign class field $H_A^+$,
normalized for the chosen sign function and with
$\Lambda_\psi=\pC\iota(A)$.
Write $\xi_\psi\in\{0,\ldots,N-1\}$ for its sign exponent, determined by
$\operatorname{LT}(\psi_a)=\iota(\operatorname{sgn}(a))^{q^{\xi_\psi}}$
for $a\in A\setminus\{0\}$, where $\operatorname{LT}$ denotes the leading
nonzero coefficient. We use the Artin convention
$\mathfrak a*\psi=\psi^{\sigma_{\mathfrak a}}$ of
\cite[Section 3.5]{HYY}. For a nonzero integral ideal $\mathfrak a$,
let $\psi_{\mathfrak a}:\psi\to\mathfrak a*\psi$ be its monic
ideal isogeny and set
\[
 d_{\mathfrak a}:=\partial\psi_{\mathfrak a},\qquad
 u_{\mathfrak a,i}
 =p_i^*\left(\frac{\psi_{\mathfrak a}(\omega_\psi)}{\omega_\psi}\right).
\]
The quotients are rational functions independent of the nonzero
special function $\omega_\psi$; see \cite[Section 2.1]{HYY}.
They are the shtuka multipliers of \cite[Section 3.1]{ANDTR17},
whose Lemma 3.1 gives their ideal-isogeny cocycle.

\begin{definition}\label{def:ideal-L-series}
For every integer $n\geqslant1$, the
\emph{multicharacter ideal-type $L$-value} $\mathbb L_{\psi}^s(n)$ is defined by
\begin{equation}\label{eq:mc-ideal-definition}
 \mathbb L_{\psi}^s(n):=\sum_{0\ne\mathfrak a\subset A}
                \frac{\prod_{i=1}^s u_{\mathfrak a,i}}
                     {d_{\mathfrak a}^{\,n}}.
\end{equation}
\end{definition}

With the compatible shtuka normalization, $\mathbb L_{\psi}^s(1)$ is the
value at $1$ of the group-ring operator $\mathcal L_s$ introduced in
\cite[Section 5]{ANDTR17}: each Artin symbol in that operator acts
trivially on $1$.

\begin{proposition}\label{prop:mc-class-decomposition}
\textup{(i)} The series \eqref{eq:mc-ideal-definition} converges for
every $n\geqslant1$. For integral representatives $\mathfrak b_j$ of the
narrow ideal class group $\operatorname{Cl}^+(A)$,
\begin{equation}\label{eq:mc-narrow}
 \mathbb L_{\psi}^s(n)=\sum_j
      \frac{\prod_i u_{\mathfrak b_j,i}}{d_{\mathfrak b_j}^{\,n}}
               L^s((\mathfrak b_j^{-1})_+;n).
\end{equation}
\par\medskip\noindent\textup{(ii)} Let
$\mathbb L_{\psi,[\mathfrak a]}^s(n)$ denote the subseries of
\eqref{eq:mc-ideal-definition} over ideals in the ordinary class
$[\mathfrak a]$. If $n\equiv s\pmod{q^N-1}$ and $\mathfrak a$ is integral, then
\begin{equation}\label{eq:mc-ordinary-class}
 \mathbb L_{\psi,[\mathfrak a]}^{s}(n)
 =-\frac{\prod_i u_{\mathfrak a,i}}{d_{\mathfrak a}^{\,n}}
               L^s(\mathfrak a^{-1};n).
\end{equation}

\end{proposition}
\begin{proof}

\emph{(i).} Write $\deg\mathfrak b:=\dim_{\Fq}(A/\mathfrak b)$.
For an integral ideal $\mathfrak b$ and $a\in A\setminus\{0\}$, put
$\epsilon_{\mathfrak b}(a)
=\iota(\operatorname{sgn}(a))^{-q^{\xi_\psi+\deg\mathfrak b}}$.
The composite $\psi_{\mathfrak b}\circ\psi_a$ has the same
kernel as $\psi_{a\mathfrak b}$ and leading coefficient
$\epsilon_{\mathfrak b}(a)^{-1}$. Monic normalization therefore gives
\[
 \psi_{a\mathfrak b}
 =\epsilon_{\mathfrak b}(a)\,\psi_{\mathfrak b}\circ\psi_a.
\]
Differentiation and evaluation at $\omega_\psi$ give, respectively,
\[
 d_{a\mathfrak b}=\epsilon_{\mathfrak b}(a)\iota(a)d_{\mathfrak b},
 \qquad
 u_{a\mathfrak b,i}=\epsilon_{\mathfrak b}(a)\chi_i(a)u_{\mathfrak b,i}.
\]
The second identity uses $\psi_a(\omega_\psi)=\chi(a)\omega_\psi$
and the fact that coefficient Frobenius fixes $\chi(a)$.

Now write $c:=a/b\in K^\times$ with $a,b\in A\setminus\{0\}$ and
$c\mathfrak b$ integral. Apply the two integral scaling laws to
$b(c\mathfrak b)=a\mathfrak b$ and divide. Since
$\deg(c\mathfrak b)-\deg\mathfrak b=-Nv_\infty(c)$, where $v_\infty$
is the normalized valuation at $\infty$, the Frobenius powers in the
sign factors agree on $\mathbb F_\infty$. Their quotient is
\[
 \alpha:=\frac{\epsilon_{\mathfrak b}(a)}{\epsilon_{c\mathfrak b}(b)}
 =\iota(\operatorname{sgn}(c))^{-q^{\xi_\psi+\deg\mathfrak b}}.
\]
The two scaling laws therefore become
\[
 d_{c\mathfrak b}=\alpha\iota(c)d_{\mathfrak b},\qquad
 u_{c\mathfrak b,i}=\alpha\chi_i(c)u_{\mathfrak b,i}.
\]
These are the principal-rescaling identities of
\cite[Lemma 2.7 and equation (2.19)]{HYY}.
Multiplying the $s$ numerator identities and dividing by the $n$th
power of the derivative identity gives
\begin{equation}\label{eq:mc-principal-rescaling}
 \frac{\prod_i u_{c\mathfrak b,i}}{d_{c\mathfrak b}^{\,n}}
 =\frac{\prod_i u_{\mathfrak b,i}}{d_{\mathfrak b}^{\,n}}
    \iota(\operatorname{sgn}(c))^{(n-s)q^{\xi_\psi+\deg\mathfrak b}}
    \frac{\prod_i\chi_i(c)}{\iota(c)^n}.
\end{equation}

Each ideal in the narrow class of $\mathfrak b_j$ has a unique
expression $c\mathfrak b_j$ with $c\in\mathfrak b_j^{-1}$ of sign one.
Summing \eqref{eq:mc-principal-rescaling} gives \eqref{eq:mc-narrow}.
There are finitely many narrow classes, so the convergence of the
sign-one element sums proves unconditional convergence of the ideal sum.

\emph{(ii).} When $(q^N-1)\mid(n-s)$, the sign factor in
\eqref{eq:mc-principal-rescaling} is one. Integral ideals in the
ordinary class of $\mathfrak a$ correspond to
$(\mathfrak a^{-1}\setminus\{0\})/\Fq^\times$.
The summand is constant on each such orbit, and the all-element sum
counts each ideal $q-1=-1$ times. This proves
\eqref{eq:mc-ordinary-class}.
\end{proof}

\begin{remark}
For $s=1$, \eqref{eq:mc-narrow} is
\cite[Corollary 2.10]{HYY}.
\end{remark}

\subsection{The Hilbert-class generating trace formula}\label{sec:trace-proof}

\begin{samepage}
\begin{theorem}\label{thm:mc-traceE}
Let $\psi$ be a sign-normalized Hayes module with
$\Lambda_\psi=\pC\iota(A)$.
\begin{enumerate}[label=(\roman*)]
\item For every $n\geqslant1$ satisfying $n\equiv s\pmod{q^N-1}$,
$C_{\psi,n}^{s}$ has rational coefficients over $H_A$ and
\begin{equation}\label{eq:mc-weight-trace}
 \mathbb L_{\psi}^s(n)\frac{\mathscr W_{\psi,s}}{\pC^{\,n}}
       =-\Tr_{H_A/K}(C_{\psi,n}^{s}).
\end{equation}
\item Consequently \eqref{eq:intro-mc-trace} holds. More precisely, if
$\mathfrak a_1,\ldots,\mathfrak a_h$ are integral representatives of
the ordinary ideal classes and $\varphi_j=\mathfrak a_j*\psi$, then
\begin{equation}\label{eq:mc-generating-class-sum}
       \sum_{j=1}^h(\mathcal P_{q^N}^{s}\mathscr C_{\varphi_j}^s)(\pC z)
 =\frac{\mathscr W_{\psi,s}}{\pC}
 \sum_{\substack{n\geqslant1\\n\equiv s\ (q^N-1)}}\mathbb L_{\psi}^s(n)z^{n-1}
\end{equation}
for $|z|<|\pC|^{-1}\min_j r_{\varphi_j}$, where
$r_{\varphi_j}=\min_{0\ne\lambda\in\Lambda_{\varphi_j}}|\lambda|$.
The finite sum on the left extends meromorphically in $z$.
\end{enumerate}
\end{theorem}
\end{samepage}
\begin{proof}

\emph{(i).} The period product $\mathscr W_{\psi,s}$ is
nonzero: for any $b\in A\setminus\Fq$, residue duality gives
\[
 \langle1/b,\Hcal_\psi(\pC)\rangle_\infty
       =\exp_\psi(\pC/\iota(b))\ne0,
\]
because $1/b\notin A$ and $\Lambda_\psi=\pC\iota(A)$.
Hence $\mathscr W_{\psi,s}$ has a nonzero scalar trivialization
in the Gauss field.

For the representatives $\mathfrak a_j$ and modules $\varphi_j$
in part~(ii), the period-isogeny identities
\cite[equations (2.4) and (3.3)]{HYY} give
\[
 \Lambda_{\varphi_j}=d_{\mathfrak a_j}\pC\iota(\mathfrak a_j^{-1}),
 \qquad
 \Hcal_{\varphi_j,i}(d_{\mathfrak a_j}\pC)
       =u_{\mathfrak a_j,i}\Hcal_{\psi,i}(\pC).
\]
Apply \eqref{eq:fractional-period-linearity} to $\varphi_j$ with
$\pi=d_{\mathfrak a_j}\pC$ and $\mathfrak j=\mathfrak a_j^{-1}$.
Combining it with the isogeny identity gives, for every
$c\in\mathfrak a_j^{-1}$,
\[
 \Hcal_{\varphi_j,i}(d_{\mathfrak a_j}\pC\iota(c))
 =\chi_i(c)u_{\mathfrak a_j,i}\Hcal_{\psi,i}(\pC).
\]
Substitute the lattice parametrization
$\lambda=d_{\mathfrak a_j}\pC\iota(c)$, with
$0\ne c\in\mathfrak a_j^{-1}$, into
Definition~\ref{def:mc-pellarin}. Taking external tensor products in the
preceding period identity then gives
\begin{align}
 C_{\varphi_j,n}^{s}
 &=\sum_{0\ne c\in\mathfrak a_j^{-1}}
   \frac{\prodext_{i=1}^s
     \Hcal_{\varphi_j,i}(d_{\mathfrak a_j}\pC\iota(c))}
        {(d_{\mathfrak a_j}\pC\iota(c))^n}\notag\\
 &=\frac{\prod_i u_{\mathfrak a_j,i}}{(d_{\mathfrak a_j}\pC)^n}
   \mathscr W_{\psi,s}
   \sum_{0\ne c\in\mathfrak a_j^{-1}}
   \frac{\prod_i\chi_i(c)}{\iota(c)^n}\notag\\
 &=\frac{\prod_i u_{\mathfrak a_j,i}}{(d_{\mathfrak a_j}\pC)^n}
   \mathscr W_{\psi,s}L^s(\mathfrak a_j^{-1};n)\notag\\
 &=-\mathbb L_{\psi,[\mathfrak a_j]}^{s}(n)
   \frac{\mathscr W_{\psi,s}}{\pC^n}.
 \label{eq:mc-class-period}
\end{align}
The last equality is \eqref{eq:mc-ordinary-class}.

By Proposition~\ref{prop:field-covariance}, $C_{\psi,n}^{s}$ has
coefficients in $H_A^+$.
Let $G_0:=\operatorname{Gal}(H_A^+/H_A)
\simeq\mathbb F_\infty^\times/\Fq^\times$.
A principal Artin symbol acts on $\psi$ by the scalar conjugation
$\psi^{\sigma_{(a)}}=z_a\psi z_a^{-1}$, where
$z_a=\iota(\operatorname{sgn}(a))^{-q^{\xi_\psi}}$
\cite[equation (3.26)]{HYY}. Hence
\[
 \sigma_{(a)}(C_{\psi,n}^{s})
       =z_a^{s-n}C_{\psi,n}^{s}=C_{\psi,n}^{s}.
\]
The congruence $(q^N-1)\mid(n-s)$ gives $z_a^{s-n}=1$.
The principal Artin symbols generate $G_0$, so $C_{\psi,n}^{s}$ is
fixed by $G_0$ and has coefficients in $H_A$.

Covariance gives
$C_{\varphi_j,n}^{s}=\sigma_{\mathfrak a_j}(C_{\psi,n}^{s})$.
The restrictions of $\sigma_{\mathfrak a_j}$ run through
$\operatorname{Gal}(H_A/K)$, and therefore
\[
 \sum_{j=1}^h C_{\varphi_j,n}^{s}
       =\Tr_{H_A/K}(C_{\psi,n}^{s}).
\]
Summing \eqref{eq:mc-class-period} gives
\eqref{eq:mc-weight-trace}.

\emph{(ii).} Substitute \eqref{eq:mc-projection} for every $\varphi_j$.
For each retained weight, \eqref{eq:mc-class-period} gives
\[
 \frac{\pC}{\mathscr W_{\psi,s}}
       \bigl(-C_{\varphi_j,n}^{s}\bigr)(\pC z)^{n-1}
       =\mathbb L_{\psi,[\mathfrak a_j]}^{s}(n)z^{n-1}.
\]
Normal convergence on every smaller closed disk permits summation over
$n$ and the finitely many classes. This proves
\eqref{eq:mc-generating-class-sum} and \eqref{eq:intro-mc-trace}.
Each $\mathscr C_{\varphi_j}^s$ is meromorphic by
Theorem~\ref{thm:main}. The projector is a finite linear combination
of argument rescalings, so the left side of
\eqref{eq:mc-generating-class-sum} is meromorphic as well.
\end{proof}

\begin{proof}[Proof of Theorem~\ref{thm:introtrace}]
Proposition~\ref{prop:L-expansion} gives the Taylor expansion and the congruence
vanishing, while Proposition~\ref{prop:field-covariance} gives rationality
over the field of definition of $\varphi$.
\end{proof}

\begin{proof}[Proof of Theorem~\ref{thm:introidealtrace}]

\emph{(i).} Proposition~\ref{prop:mc-element-normalization}(i)
gives the coefficient normalization. Theorem~\ref{thm:introtrace} then gives
\eqref{eq:intro-element-generating} after the substitution $w=\pC z$.

\emph{(ii).} Theorem~\ref{thm:mc-traceE}(ii) gives the finite ordinary-class generating
sum, including \eqref{eq:intro-mc-trace}, and its meromorphic continuation.
\end{proof}

\Needspace{7\baselineskip}
\section{Examples and parameter specializations}\label{sec:curve-examples}

\subsection{Low-degree examples}

\begin{example}\label{ex:small-m}
For $m=|S|<2q$, Proposition~\ref{prop:trichotomy} gives
$d_{\max}\leqslant0$, so every nonzero $G_S$ is constant. If $G_S\ne0$,
Lemma~\ref{lem:support} at $d=0$ gives $m=q+(q-1)k$ with $k\geqslant0$.
Since $m<2q$,
\[
 (q-1)k<q,\qquad\text{so }k<\frac q{q-1}\leqslant2.
\]
Thus $k\in\{0,1\}$, so only $m=q$ and $m=2q-1$ can give nonzero corrections;
all other $m<2q$ give $G_S=0$. In these two cases,
Corollary~\ref{cor:degrees} gives \eqref{eq:leading} and
\eqref{eq:second-leading}.
\end{example}

\begin{example}\label{ex:first-nonconstant}
For $m=|S|=2q$, Proposition~\ref{prop:trichotomy} gives
$\deg_XG_S=d_{\max}=1$.
If $q>2$, Lemma~\ref{lem:support} at $d=0$ would require
$2q=q+(q-1)k$, hence $k=q/(q-1)\notin\mathbb Z$. Thus \eqref{eq:leading}
with $d_{\max}=1$ gives the first nonconstant correction
\[
 \mathfrak G_S(X)=X\prodext_{i\in S}\Delta_{\varphi,i}
 \qquad(q>2).
\]

If $q=2$, then $|S|=4$, and \eqref{eq:leading} with $d_{\max}=1$
gives the same degree-one coefficient. For the constant term,
Lemma~\ref{lem:support} at $d=0$ gives $4=2+k$, hence $k=2$.
The depth bound in Lemma~\ref{lem:chain-contribution}(ii)
therefore gives $n_i\leqslant2$.
Only $x_0,x_1,x_2$ can be nonzero. The two equations in
Lemma~\ref{lem:depth-profile} become $x_0+x_1+x_2=4$ and
$x_0+x_1/2+x_2/4=2$, so $2x_1+3x_2=8$. The profiles are $(0,4,0)$ and
$(1,1,2)$; up to permutation, their depth tuples are $(1,1,1,1)$ and
$(0,1,2,2)$, respectively. Consequently,
\eqref{eq:coefficient-formula} gives the full polynomial
\[
 \begin{aligned}
 \mathfrak G_S(X)
 &=X\prodext_{i\in S}\Delta_{\varphi,i}
   +\prodext_{i\in S}\mathfrak{B}_{i}^{[1]}\\
 &\quad+\sum_{\substack{i,j\in S\\i\ne j}}
   \Delta_{\varphi,i}\boxtimes \mathfrak{B}_{j}^{[1]}\boxtimes
   \prodext_{\ell\in S\setminus\{i,j\}}\mathfrak{B}_{\ell}^{[2]}
 \qquad(q=2,\ |S|=4),
 \end{aligned}
\]
where \eqref{eq:B-definition} gives
$\mathfrak{B}_{\ell}^{[2]}=(f_\ell \mathfrak{B}_{\ell}^{[1]})^{(-1)}$. The last sum
has $12$ terms, so the constant coefficient has $1+12=13$ terms in the
differential basis.
\end{example}

\subsection{The Carlitz correction polynomials}

\begin{example}\label{ex:carlitz-q3}

For ${\mathcal X}=\mathbb P^1$ and the Carlitz module over
$A:=\Fq[\theta]$, denote by $H_S$ the correction polynomial indexed
by $\Sigma_s\setminus S$ in \cite[Theorem 2]{PP16}.
For $|t_i|\leqslant1$ ($i\in S$), these polynomials satisfy the canonical
sum formula of \cite[Section 4.1]{PP16},
\[
 H_S(e_C(z))=\pC^{-1}\sum_{a\in\Fq[\theta]}
       \frac{\prod_{i\in S}f_{t_i}(a-z)}{z-a}.
\]

Using Notation~\ref{not:carlitz-normalization} and substituting
$w=\pC z$ and $\lambda=\pC a$ into
Proposition~\ref{prop:lattice-limit}(ii) gives
\[
 G_S(e_C(z))=(-1)^{|S|}\pC^{-1}
 \sum_{a\in\Fq[\theta]}\frac{\prod_{i\in S}f_{t_i}(a-z)}{z-a}
 =(-1)^{|S|}H_S(e_C(z)).
\]
Hence $G_S=(-1)^{|S|}H_S$ and, since $\nu_S:=\prodext_{i\in S}dt_i$,
$\mathfrak{G}_S=(-1)^{|S|}H_S\prodext_{i\in S}dt_i$.

Proposition~\ref{prop:trichotomy}(ii)--(iii) therefore shows that
$H_S$ is constant exactly when $|S|\notin\Nq$, and
Corollary~\ref{cor:simultaneous-constancy} shows that all $H_S$ for
nonempty $S\subseteq\Sigma_s$ are constant exactly when $s<2q$.

At the first nonconstant size $|S|=2q$,
Example~\ref{ex:first-nonconstant} gives $H_S(X)=X$ for $q>2$ and
$H_S(X)=X+g_{S,0}$ for $q=2$. In the latter case, let $\sigma_j$ be the
$j$th elementary symmetric polynomial in the four parameters. Then
\[
 g_{S,0}=\sigma_2\sigma_3+\sigma_1\sigma_4+\sigma_4
             +\theta(\sigma_2+\sigma_3)+\theta^2.
\]
To obtain this expression, substitute
$\mathfrak{B}^{[1]}(t)=(t-\theta^{1/2})\,dt$ and
$\mathfrak{B}^{[2]}(t)=(t-\theta^{1/2})(t-\theta^{1/4})\,dt$ into
the thirteen terms of Example~\ref{ex:first-nonconstant}.
After expansion, the terms with exponent pattern $(2,2,1,0)$ sum to
$\sigma_2\sigma_3+\sigma_1\sigma_4$ in characteristic two; the remaining
terms are $\sigma_4+\theta(\sigma_2+\sigma_3)+\theta^2$.
\end{example}

\subsection{Parameter recursions and diagonal corrections}
\label{sec:carlitz-carries}

Use the Carlitz normalization of
Example~\ref{ex:carlitz-q3}. To record parameter specializations, write
\[
 H(t_1,\ldots,t_m;X):=H_{\Sigma_m}(X)\qquad(m\geqslant1),
 \qquad H(\,;X):=H_\varnothing(X):=X^{-1}.
\]
Repeated arguments mean that the corresponding geometric parameters are
set equal; for example, $H(t,t,t;X)$ is the full diagonal specialization
of $H_{\Sigma_3}$.
In this normalization,
$\mathfrak{B}^{[n]}(t)=\prod_{j=1}^n(t-\theta^{q^{-j}})\,dt$.
Theorem~\ref{thm:formula} therefore gives scalar coefficients polynomial
in the geometric parameters, which can be specialized at any values.
The identity $a(t)^q=a(t^q)$ for $a\in\Fq[\theta]$ combines $q$ equal
evaluation characters into one. For the correction polynomials, the
Carlitz exponential equation gives the additional term
$-XH(t_1,\ldots,t_m;X)$ in \eqref{eq:carlitz-carry}.

\Needspace{10\baselineskip}
\begin{proposition}\label{prop:carlitz-carries}
\leavevmode 
\begin{enumerate}[label=(\roman*)]
\item For every $m\geqslant0$ and all parameters
$t_1,\ldots,t_m,t$,
\begin{equation}\label{eq:carlitz-carry}
\begin{aligned}
 &H(t_1,\ldots,t_m,t,\ldots,t;X)\\
 &\quad=(t^q-\theta)H(t_1,\ldots,t_m,t^q;X) -XH(t_1,\ldots,t_m;X),
\end{aligned}
\end{equation}
where $t$ occurs $q$ times after $t_m$.
\item  For every $r\geqslant1$,
\begin{equation}\label{eq:carlitz-diagonal}
  H(\underbrace{t,\ldots,t}_{q^r};X)
   =-\sum_{j=0}^{r-1}
       \left(\prod_{h=j+1}^{r-1}(t^{q^r}-\theta^{q^h})\right)
       X^{q^j-1},
\end{equation}
where the argument $t$ appears $q^r$ times. In particular, its degree is
$q^{r-1}-1$ and its leading coefficient is $-1$, even after any further
specialization of $t$. 
\end{enumerate}
\end{proposition}

\begin{proof}

\emph{(i).} The Carlitz exponential equation and the definitions
in Notation~\ref{not:carlitz-normalization} give
$e_k^q=e_{k-1}-\theta e_k$ for $k\geqslant1$.
Consequently, reindexing the series defining $f_t$ gives
\begin{align}
 f_t(z)^q
 &=\sum_{k\geqslant1}t^{q(k-1)}
   \bigl(e_{k-1}(z)-\theta e_k(z)\bigr)\notag\\
 &=e_C(z)+t^q f_{t^q}(z)-\theta f_{t^q}(z)\notag\\
 &=(t^q-\theta)f_{t^q}(z)+e_C(z).
 \label{eq:carlitz-parameter-frob}
\end{align}

For parameters in the closed unit disk, $f_t(a)=\omega(t)a(t)$ is
uniformly bounded in $a\in A$, and $f_t(z)$ is bounded on closed
$z$-disks. Thus the canonical sums converge normally away from $A$.
Evaluation on repeated parameters is a bounded map of Tate algebras
and commutes with these sums.
Put $X=e_C(z)$. For $a\in A$, the identity $\ker e_C=A$ gives
$e_C(a-z)=-X$, so \eqref{eq:carlitz-parameter-frob} becomes
\[
 f_t(a-z)^q=(t^q-\theta)f_{t^q}(a-z)-X.
\]
Substitution into the normally convergent canonical sum gives
\begin{align*}
 H(t_1,\ldots,t_m,\underbrace{t,\ldots,t}_{q};X)
 &=\pC^{-1}\sum_{a\in A}
 \frac{\bigl(\prod_{i=1}^m f_{t_i}(a-z)\bigr)f_t(a-z)^q}{z-a}\\
 &=(t^q-\theta)\pC^{-1}\sum_{a\in A}
 \frac{\bigl(\prod_{i=1}^m f_{t_i}(a-z)\bigr)f_{t^q}(a-z)}{z-a}\\
 &\quad-X\pC^{-1}\sum_{a\in A}
 \frac{\prod_{i=1}^m f_{t_i}(a-z)}{z-a}\\
 &=(t^q-\theta)H(t_1,\ldots,t_m,t^q;X)-XH(t_1,\ldots,t_m;X).
\end{align*}
Here the last line uses the canonical-sum formula in
Example~\ref{ex:carlitz-q3}. For $m=0$, the scalar partial fraction
identity supplies the empty-product sum $H_\varnothing(X)=X^{-1}$,
so the same calculation applies. Both sides of
\eqref{eq:carlitz-carry} are polynomials in $X$, since
$XH_\varnothing(X)=1$.
They agree at $X=e_C(z)$ for $z\notin A$ and parameters in the closed
unit disk. The infinite image of $e_C$ gives equality as polynomials
in $X$. Since their coefficients are polynomial in the parameters,
equality on Cartesian powers of the infinite unit disk extends the
identity to all parameter values.

\emph{(ii).} Induction on $r$ gives
\begin{align*}
 f_t(z)^{q^r}
 &=\left(\prod_{h=0}^{r-1}(t^{q^r}-\theta^{q^h})\right)
        f_{t^{q^r}}(z)\\
 &\quad+\sum_{j=0}^{r-1}
       \left(\prod_{h=j+1}^{r-1}(t^{q^r}-\theta^{q^h})\right)
       e_C(z)^{q^j}.
\end{align*}
For $r=1$ this is \eqref{eq:carlitz-parameter-frob}.
Taking the $q$th power shifts the indices $h$ and $j$ by one;
applying \eqref{eq:carlitz-parameter-frob} to $f_{t^{q^r}}$
supplies the new factor with $h=0$ and the new summand with $j=0$.

Apply the canonical sum with all $q^r$ parameters equal to $t$.
The first term contributes zero because the one-factor correction
vanishes. For each $j$, the term involving $e_C^{q^j}$ contributes its
displayed coefficient times
\[
 \pC^{-1}\sum_{a\in A}\frac{e_C(a-z)^{q^j}}{z-a}
   =(-X)^{q^j}X^{-1}=-X^{q^j-1}.
\]
Substitution gives
\eqref{eq:carlitz-diagonal} first for $|t|\leqslant1$ and then, by
polynomiality, for every $t$. For $j=r-1$, the product is empty, so
the leading coefficient is $-1$.
\end{proof}

\subsection{A nonprincipal lattice}
Take $q=3$, $\rho(t):=t^3-t-1$, and
$A:=\Gamma(\mathbb P^1\setminus P_\rho,\mathcal O)$,
where $P_\rho$ is the closed point defined by $\rho$.
The polynomial $\rho$ is irreducible, since it has no root in $\mathbb F_3$.
The divisor class sequence identifies $\operatorname{Cl}(A)$ with
$\operatorname{Pic}(\mathbb P^1)/\langle[P_\rho]\rangle
\simeq\mathbb Z/3\mathbb Z$, since $\deg P_\rho=3$.
It is generated by the class of
$I_{\infty_0}=(1/\rho,t/\rho,t^2/\rho)$, the ideal of the usual rational
point $\infty_0$.

The projective-line module $\varphi$ in
\cite[Section 5.1 and Lemma 5.2]{HYY} has lattice
$\kappa\iota(I_{\infty_0}^2)$ for a nonzero scalar $\kappa\in\Cinf$
and shtuka function
\[
 f(t)=\frac1{t-\eta}-\frac1{\theta-\eta},
 \qquad \rho(\eta)=0,\quad \theta:=\iota(t).
\]
Since $I_{\infty_0}^2$ is nonprincipal, this period lattice
is nonfree. Choose the root $\eta$ compatible with the infinite-place
embedding, so that $|\theta-\eta|<1$. Write
$\eta_{-1}=\eta^{3^2}$, $\eta_{-2}=\eta^3$, and $\beta:=\theta^{1/3}$.
\par\Needspace{7\baselineskip}
The Cauchy differential and normalized shtuka differential of $\varphi$ are
\begin{align*}
 C_\varphi&=\left(\frac1{t-\theta}-\frac1{t-\eta_{-1}}\right)dt,\\
 \Delta_\varphi&=-\frac{\beta-\eta_{-2}}
 { (\beta-\eta_{-1})(t-\eta_{-1})(t-\eta_{-2})}\,dt.
\end{align*}
A differential on $\mathbb P^1$ with only simple poles at $t=\theta$
and $t=\eta_{-1}$ is determined by its residue at $t=\theta$.
Multiplying $C_\varphi$ by
$f$ and taking one inverse coefficient twist gives $\Delta_\varphi$.

On six independent geometric factors, $m=6=2q$, so
Example~\ref{ex:first-nonconstant} gives the first nonconstant correction
\[
 \mathfrak{G}_{\{1,\ldots,6\}}(X)=X\prodext_{i=1}^6\Delta_{\varphi,i}.
\]

\subsection{A constant correction beyond the simultaneous-constancy threshold}
Take $q=3$ and a factor set $S$ with $m=|S|=7$.
The degree formula gives
\[
 d_{\max}=7-1-2\left\lceil\frac73\right\rceil=0,
\]
so Proposition~\ref{prop:trichotomy} gives
$\mathfrak G_S(X)=\mathfrak g_{S,0}\ne0$, even though $m>2q$.

To compute this constant from \eqref{eq:coefficient-formula}, we determine
$\mathbb E_{S,0}$. The depth bound in Lemma~\ref{lem:chain-contribution}(ii)
gives $n_i\leqslant(7-3)/(3-1)=2$, so only $x_0,x_1,x_2$ can be nonzero.
Lemma~\ref{lem:depth-profile} gives
\[
 x_0+x_1+x_2=|S|=7,\qquad
 x_0+\frac{x_1}{3}+\frac{x_2}{9}=q(0+1)=3.
\]
Thus $x_1=10-4x_0$ and $x_2=3x_0-3$. Their nonnegativity forces
$x_0\in\{1,2\}$, giving the two depth profiles that exhaust
$\mathbb E_{S,0}$:
\[
 (x_0,x_1,x_2)=(1,6,0)\quad\mbox{or}\quad(2,2,3).
\]
The corresponding depth tuples are, up to permutation,
\[
      (0,1,1,1,1,1,1)\quad\mbox{or}\quad(0,0,1,1,2,2,2).
\]
Substitution in \eqref{eq:coefficient-formula} therefore gives
\[
\begin{aligned}
\mathfrak G_S(X)=\mathfrak g_{S,0}
={}&\sum_{(n_i)\in\mathbb E_{S,0}}
      \prodext_{i\in S}\mathfrak{B}_i^{[n_i]}\\
={}&\sum_{a\in S}\mathfrak{B}_a^{[0]}\boxtimes
      \prodext_{i\in S\setminus\{a\}}\mathfrak{B}_i^{[1]}\\
&+\sum_{\substack{S_0,S_1\subseteq S\\
                  S_0\cap S_1=\varnothing\\
                  |S_0|=|S_1|=2}}
 \Biggl[
 \left(\prodext_{i\in S_0}\mathfrak{B}_i^{[0]}\right)
 \boxtimes
 \left(\prodext_{i\in S_1}\mathfrak{B}_i^{[1]}\right) \boxtimes
 \left(\prodext_{i\in S\setminus(S_0\cup S_1)}\mathfrak{B}_i^{[2]}\right)
 \Biggr].
\end{aligned}
\]
Here \eqref{eq:B-definition} gives $\mathfrak{B}_i^{[0]}=\Delta_{\varphi,i}$ and
$\mathfrak{B}_i^{[2]}=(f_i\mathfrak{B}_i^{[1]})^{(-1)}$.

The first sum has $7$ terms. The second has
$\binom72\binom52=210$, obtained by choosing the disjoint two-element
sets $S_0$ and $S_1$. Thus $|\mathbb E_{S,0}|=217$, and these distinct
external basis elements give a nonzero constant by
Theorem~\ref{thm:formula}(i). The seven-factor correction is therefore
constant, although the same family contains nonconstant corrections
on six factors.

\section*{Acknowledgments}
The authors used ChatGPT and OpenAI Codex (OpenAI) for manuscript
organization, English-language revision, LaTeX preparation, and
preliminary exploration and checking of mathematical arguments. The
authors reviewed the resulting text and arguments and take responsibility
for the mathematical content and references.

\end{document}